\documentclass[onecolumn, draftcls]{IEEEtran}
\usepackage{graphicx}
\graphicspath{{imgs/}}
\usepackage{float}
\usepackage{amssymb,latexsym,amsmath,mathtools}
\numberwithin{equation}{section}
\usepackage{xcolor,comment}
\usepackage[center]{caption}
\usepackage{color}
\usepackage{tikz}
\usetikzlibrary{shapes.geometric, arrows}
\tikzstyle{arrow} = [thick,->,>=stealth]
\usepackage{algpseudocodex}
\usepackage{tabularx}
\usepackage{pifont}
\newcommand{\cmark}{\ding{51}}%
\newcommand{\xmark}{\ding{55}}%

\newtheorem{theorem}{Theorem}[section]
\newtheorem{proposition}[theorem]{Proposition}
\newtheorem{corollary}[theorem]{Corollary}
\newtheorem{lemma}[theorem]{Lemma}

\newtheorem{definition}[theorem]{Definition}

\DeclareMathOperator*{\argmin}{argmin}

\allowdisplaybreaks

\begin{document}
\title{Sliding Window Codes: Near-Optimality and Q-Learning for Zero-Delay Coding

\thanks{This work was supported by the Natural Sciences and Engineering Research Council of Canada. The material in this paper has been presented in part (without the associated technical analysis) at the 2023 Conference on Decision and Control and the 2024 American Control Conference.}}

\author{Liam Cregg, Fady Alajaji, Serdar Y{\"u}ksel\\
\thanks{The first author is with the Department of Information Technology and Electrical Engineering, ETH Z{\"u}rich, Switzerland. The second and third authors are with the Department of Mathematics and Statistics, Queen's University, Kingston, Canada (e-mail: lcregg@ethz.ch; fa@queensu.ca; yuksel@queensu.ca).}
}

\maketitle
\begin{abstract}%
    We study the problem of zero-delay coding for the transmission of a Markov source over a noisy channel with feedback and present a reinforcement learning solution which is guaranteed to achieve near-optimality. To this end, we formulate the problem as a Markov decision process (MDP) where the state is a probability-measure valued predictor/belief and the actions are quantizer maps. This MDP formulation has been used to show the optimality of certain classes of encoder policies in prior work, but their computation is prohibitively complex due to the uncountable nature of the constructed state space and the lack of minorization or strong ergodicity results. These challenges invite rigorous reinforcement learning methods, which entail several open questions: can we approximate this MDP with a finite-state one with some performance guarantee? Can we ensure convergence of a reinforcement learning algorithm for this approximate MDP? What regularity assumptions are required for the above to hold? We address these questions as follows: we present an approximation of the belief MDP using a sliding finite window of channel outputs and quantizers. Under an appropriate notion of predictor stability, we show that policies based on this finite window are near-optimal, in the sense that the lowest distortion achievable by such a policy approaches the true lowest distortion as the window length increases. We give sufficient conditions for predictor stability to hold. Finally, we propose a Q-learning algorithm which provably converges to a near-optimal policy and provide a detailed comparison of~the sliding finite window scheme with another approximation scheme which quantizes the belief MDP in a nearest neighbor fashion.
\end{abstract}

\begin{IEEEkeywords}
    Zero-delay coding, reinforcement learning, Q-learning, stochastic control, Markov decision processes, quantization, source-channel coding, noisy channels with feedback
\end{IEEEkeywords}

\section{Introduction}\label{sectionIntro}
The zero-delay coding problem involves compressing and transmitting an information source at a fixed rate over a noisy channel with feedback and without delay, while minimizing the expected distortion at the receiver. This zero-delay restriction is of practical relevance in many applications, including live-streaming \cite{badr2013streaming,rudow2022streaming,khisti2014streaming,khisti16}, real-time tracking and estimating processes over erasure links~\cite{khina2018tracking} and real-time sensor networks \cite{Gastpar,GastparMAC}. However, this restriction means that classical Shannon-theoretic methods~\cite{Shannon}, which require collecting large sequences of source symbols and compressing them at once, are not viable as they induce a large delay. From an information-theoretic perspective, several strategies have been used to approach this problem, including mutual information constraints, entropy coding, and Shannon lower bounding techniques. Studies to this end include \cite{Silva1, silva2015characterization, stavrou2018zero}. Within the context of linear systems, \cite{banbas89,TatikondaSahaiMitter, tanaka2016semidefinite,
    derpich2012improved,stavrou2021asymptotic} use sequential rate-distortion theory. Some of these works give applicable codes for zero-delay coding for Gaussian sources over additive-noise Gaussian channels, and some give upper and/or lower performance bounds; see \cite{stavroutanaka2019time} and
\cite{kostina2019rate} for further studies, and see~\cite{Chakravorty2017, Guo2021, Soleymani2024} for studies in the context of channels with delay and other communication constraints.

Furthermore, learning theoretic methods have attracted significant interest in source-channel coding theory both in the classical literature and the recent literature, see for example \cite{Pollard,linder1994rates,linder2002learning} for the noiseless channel (quantization) case among several classical results, although usually restricted to independent and identically distributed (i.i.d.) sources. We note that our results are directly applicable for i.i.d.\ sources as well, since an optimal zero-delay code for an i.i.d.\ source is a memoryless code \cite{WalrandVaraiya,Witsenhausen,wood2016optimal} (see also, for related discussions in a different causal coding context \cite{NeuhoffGilbert,GaarderSlepian,Weissman,AsnaniWeissman}).
More recently, deep learning is employed to construct powerful joint source-channel codes (see~\cite{gunduz2019machine,Farsad,Kurka}), and reinforcement learning is used as a tool to estimate feedback capacity in~\cite{PermuterFiniteState,PermuterIsing}. 
Although effective in practice, these machine learning methods are generally experimental and do not provide a formal proof of convergence or optimality. Conversely, our reinforcement learning approach will be rigorously shown to converge to \emph{near-optimality}, in the sense that, for any \(\epsilon > 0\), there exists some window size $N$ such that the resulting policy achieves distortion no further than \(\epsilon\) from the optimal distortion (a formal definition is given in Definition \ref{definition:near-optimal}).

There have been several studies about the zero-delay coding problem using stochastic control techniques. In particular,~\cite{Witsenhausen,WalrandVaraiya,Teneketzis} consider Markov sources with finite alphabets and finite time horizons and show optimality of structured classes of policies. Similar optimality and existence results are presented for infinite time horizons in~\cite{wood2016optimal} (with feedback) and~\cite{MahTen09} (without feedback). The continuous-alphabet infinite-horizon case is examined in~\cite{ghomi2021zero}, although only over a noiseless channel. These results often rely on formulating the problem as a Markov decision process (MDP) in order to utilize existing results from stochastic control theory, such as dynamic programming and value iteration methods (see \cite{Hernandez, yuksel2023control} for detailed information on such methods). However, in the formulation of the MDP, these results utilize a state space that is probability measure-valued (this state is often called the ``predictor'' in the literature) and an action space involving quantizers. These spaces are computationally difficult to work with, both in terms of complexity and implementation. Thus, while numerous existence and structural results have been established for this problem, the explicit development of effective coding schemes for a given zero-delay coding problem is still an open problem. 

In this paper, we present an approximation method to simplify the resulting MDP, and we use this approximation to obtain near-optimal coding schemes for the zero-delay coding problem via a reinforcement learning approach. We emphasize that we provide guaranteed approximation and convergence results. In particular, we build on methods used in \cite{ kara2021convergence}, which were originally used to study partially observed Markov decision processes (POMDPs). Based on these methods, we introduce a practical sliding finite window method and present several mathematical and algorithmic results on its near-optimal performance. We also compare this method with the approach used in \cite{creggZeroDelayNoiseless}, which uses a nearest neighbor quantization of the probability measure-valued state space. Although \cite{creggZeroDelayNoiseless} only studied the noiseless channel (quantization) version of the problem, extensions to the noisy channel case follow with little additional effort as we discuss in the paper.

We emphasize that our new sliding finite window method has several advantages over the belief quantization scheme. In particular, it is easier to implement, computationally less complex, and valid for any initialization. This comes at a cost of some additional Dobrushin coefficient conditions on the source and channel, but we note that the conditions provided here are sufficient, not necessary. These differences can be seen by comparing Theorems~\ref{theorem:Qlearning2} and~\ref{theorem:Qlearning1}, and are further detailed in Section~\ref{section:comparison}.

The sliding finite window code that we propose can be viewed as a finite state code \cite{pollara1988finite, KiefferDunham} where the states are the realizations of a finite window of encoder maps and outputs. We are not aware of such a construction in the literature with rigorous performance guarantees in the context of optimal real-time coding and in particular on an associated learning theoretic study. We also note the similarities to trellis or convolutional codes \cite{Ungerboeck, Calderbank, Johannesson}; our method can be in some sense be seen as a zero-delay analog of these codes, but our analysis is based on stochastic control and thus entirely different than that of trellis codes.

\subsection*{Main Contributions}
\begin{itemize}
    \item To our knowledge there is no prior study on near-optimality of sliding window coding schemes for the zero-delay coding problem with an explicit rate of convergence (involving all window lengths).
    To this end, in Theorem \ref{theorem:finite window_near_optimal}, we provide an approximation result for the zero-delay coding problem which bounds the sub-optimality gap of sliding finite window policies by a predictor stability term.

    \item Under Dobrushin coefficient conditions, we then show that these sliding finite window policies become near-optimal for the zero-delay coding problem in Corollary \ref{corollary:dobrushin_finite_window}, with an explicit performance bound. We note that these conditions hold when the channel is ``noisy enough'' or when the source has some mixing conditions. To our knowledge, this gives the first rigorous proof of optimality of sliding finite window policies for the zero-delay coding problem.

    \item We provide a (reinforcement) Q-learning algorithm which allows for the computation of these sliding finite window policies and rigorously show its convergence in Theorem \ref{theorem:Qlearning2}.

    \item Finally, we provide both theoretical and experimental comparisons of our sliding finite window method with a nearest neighbor approximation scheme presented in \cite{creggZeroDelayNoiseless} (and adapted and generalized to the noisy channel setting).
\end{itemize}

{\bf Notation.} In general, we will denote random variables by capital letters and their realizations by lowercase letters. There are a few exceptions to this; in particular we will always use lowercase \(\pi\) and uppercase \(Q\) in order to avoid a conflict of notation with existing results in the literature. It will be clear from the context for these variables whether we are referring to a random variable or its realization. To denote the set of probability measures over a measurable space \((\mathcal{X}, \mathcal{B}(\mathcal{X}))\), we use \(\mathcal{P}(\mathcal{X})\), and to denote a contiguous tuple of random variables \((X_0, X_1, \ldots, X_n)\) we will use the notation \(X_{[0,n]}\) (and its realization by \(x_{[0,n]}\)). Probabilities and expectations
  will be denoted by $P$ and $\mathbf{E}$, respectively. When the relevant distributions depend on some parameters, we include these in the superscript and/or subscript. For probabilities involving finite spaces, we will often use the shorthand $P(y_t|x_t) = P(Y_t = y_t | X_t = x_t)$, or simply \(P(y|x)\), when the time index is not important. Also note that, even for a finite set \(\mathcal{Y}\), we sometimes write for consistency of notation \(\sum_{\mathcal{Y}} f(y) P(y | x) = \int_\mathcal{Y} f(y) P(dy | x)\), where we use the counting measure over \(\mathcal{Y}\).

\section{Preliminaries: Optimal Coding Problem and its MDP Formulation}\label{section:preliminaries}
\subsection{Optimal Zero-Delay Coding and Existence of an Optimal Policy}
We recall a few basic definitions and results among the large body of research on Markov chains \cite{Hernandez-Lerma2003}. Let $(a_t)_{t \ge 0}$ be an $A$-valued Markov chain, where $(A,\cal{B})$ is a measurable space. Denote its transition kernel by $K(B | a) \coloneqq P(a_{t+1} \in B | a_t = a)$ for any measurable set $B \in \cal{B}$, and its $t$-fold transition kernel by 
\begin{align}
    K^t(B|a) \coloneqq \int_{A^{t-1}} K(da_1|a)K(da_2|a_1) \ldots K(da_{t-1}|a_{t-2})K(B|a_{t-1}).
\end{align}

\begin{definition}
    We say a probability measure $\mu$ is \emph{invariant} for kernel $K$ (equivalently, for the Markov chain $(a_t)_{t \ge 0}$) if
    \begin{align}
        \int_A K(B|a)\mu(da) = K(B) \text{ for all } B \in \cal{B}.
    \end{align} 
\end{definition}

\begin{definition}
    Let $\psi$ be a nontrivial $\sigma$-finite measure over $\cal{B}$. Then $(a_t)_{t \ge 0}$ is $\psi$\emph{-irreducible} if 
    \begin{align}
        \sum_{t \ge 0} K^t(B|a) > 0 \text{ for all } a \in A \text{ whenever } \psi(B) > 0.
    \end{align}
\end{definition}
When $A$ is countable, we let $\psi$ be the counting measure and simply call the Markov chain irreducible. One classical result in Markov chain theory is that an irreducible Markov chain on a finite space has a unique invariant measure~\cite{Hernandez-Lerma2003}.

Let our information source be a Markov process \((X_t)_{t \ge 0}\) taking values in \(\mathcal{X}\), which we assume is finite. Let \(T(x' | x) \coloneqq P(X_{t+1} = x' | X_t = x)\) be its transition kernel, which we assume is irreducible, and thus admits a unique invariant measure. Denote this unique invariant measure by \(\zeta\). Let \(X_0 \sim {\pi_0}\) (we also call \({\pi_0}\) the prior). Let \(\mathcal{M}\) and \(\mathcal{M'}\) be the input and output alphabets of the memoryless noisy channel, which we assume are finite, and let \((M_t)_{t \ge 0}\) and \((M'_t)_{t \ge 0}\) be the respective processes. We denote the channel kernel by \(O(m' | m) \coloneqq P(M'_t = m' | M_t = m)\) for $m' \in \mathcal{M}'$ and $m \in \mathcal{M}$. Finally, let \(\hat{\mathcal{X}}\) be some finite set of reconstruction values, and let \((\hat{X}_t)_{t \ge 0} \subset \hat{\mathcal{X}}\) be the corresponding sequence of reproductions.

Consider sequences of functions \((\gamma^e_t)_{t \ge 0}\), which we call the encoder policy, and \((\gamma^d_t)_{t \ge 0}\), which we call the decoder policy. In addition to the current source symbol, the encoder has access to all past source symbols and channel inputs, and all past channel outputs in the form of feedback. In addition to the current channel output, the decoder has access to all previous channel outputs. That is, \((\gamma^e_t)_{t \ge 0}\) and \((\gamma^d_t)_{t \ge 0}\) are such that
\begin{align}
    \gamma^e_t : \mathcal{X}^{t+1} \times \mathcal{M}^t \times \left(\mathcal{M'}\right)^t & \to \mathcal{M} & \gamma^d_t : \left(\mathcal{M'}\right)^{t+1} & \to \hat{\mathcal{X}} \\
    (X_{[0,t]}, M_{[0,t-1]}, M'_{[0,t-1]})                                                 & \mapsto M_t     & M'_{[0,t]}                                   & \mapsto \hat{X}_t.
\end{align}
Our setup can be seen in Figure \ref{fig:coding_basic}. Note that the initial distribution $\pi_0$ and the policies $(\gamma^e_t, \gamma^d_t)_{t \ge 0}$ induce a joint distribution on $(X_t, M_t, M'_t, \hat{X}_t)_{t \ge 0}$; this follows from the Ionescu-Tulcea Theorem, see e.g., \cite[Proposition C.10]{Hernandez}.

\tikzstyle{block} = [rectangle, rounded corners, minimum width=.5cm, minimum height=.5cm,align=center, draw=black]
\tikzstyle{openblock} = [minimum width=.5cm, minimum height=.5cm,align=center]

\begin{figure}[h!]
    \centering
    \begin{tikzpicture}[node distance=3cm]
    \node (source) [openblock] {$X_t$};
    \node (encoder) [block, right of=source] {Encoder};
    \node (channel) [block, right of=encoder] {Channel};
    \node (decoder) [block, right of=channel] {Decoder};
    \node (reconstruction) [openblock, right of=decoder] {$\hat{X}_t$};
    \node (m) [openblock, above of=encoder, xshift=1.5cm, yshift=-2.5cm] {$M_t$};
    \node (m') [openblock, above of=decoder, xshift=-1.5cm, yshift=-2.5cm] {$M'_t$};
    \draw [arrow] (source) -- (encoder);
    \draw [arrow] (encoder) -- (channel);
    \draw [arrow] (channel) -- (decoder);
    \draw [arrow] (decoder) -- (reconstruction);
    \draw [arrow] (channel.east) -| +(.5,-1) -| (encoder);
\end{tikzpicture}
    \caption{Source-channel coding with feedback}
    \label{fig:coding_basic}
\end{figure}

We consider two performance criteria for the zero-delay coding problem. We wish to find encoder and decoder policies such that one of the following distortion quantities is minimized: the discounted distortion,
\begin{equation}
    J_\beta(\pi_0, \gamma^e, \gamma^d) \coloneqq \mathbf{E}^{\gamma^e, \gamma^d}_{\pi_0} \left[ \sum_{t=0}^\infty \beta^t d(X_t, \hat{X}_t)\right]\label{eq:dis_cost},
\end{equation}
or the average distortion,
\begin{equation}
    J(\pi_0, \gamma^e, \gamma^d) \coloneqq \limsup_{T \to \infty} \mathbf{E}^{\gamma^e, \gamma^d}_{\pi_0} \left[ \frac{1}{T} \sum_{t=0}^{T-1} d(X_t, \hat{X}_t)\right]\label{eq:avg_cost},
\end{equation}
where \(d : \mathcal{X} \times \hat{\mathcal{X}} \to \mathbb{R}_+\) is a given distortion function and \(\beta \in (0,1)\) is a given discount factor.

We refer to the minimization of~\eqref{eq:dis_cost} as the discounted distortion problem and of~\eqref{eq:avg_cost} as the average distortion problem. Note that for a fixed encoder policy \(\gamma^e\), it is straightforward to show that the optimal decoder policy, for all \(t \ge 0\), is given by
\begin{equation}
    \hat{X}_t = \gamma^{d*}_t(M'_{[0,t]}) \coloneqq \argmin_{\hat{x} \in \hat{\mathcal{X}}}\mathbf{E}^{\gamma^e}_{\pi_0}\left[ d(X_t,\hat{x}) | M'_{[0,t]} \right].\label{eq:decoder}
\end{equation}
Accordingly, we assume that we use an optimal decoder policy for a given encoder policy. We then denote by \(\Gamma\) the set of all encoder policies, and refer to the resulting encoder-decoder policy $(\gamma^e, \gamma^{d*})$ simply as $\gamma$. Then \eqref{eq:dis_cost} and \eqref{eq:avg_cost} become $J_\beta(\pi_0, \gamma)$ and $J(\pi_0, \gamma)$, respectively. We denote the minimal distortions by 
\begin{align}
    J^*_\beta(\pi_0) & \coloneqq \inf_{\gamma} J_\beta(\pi_0, \gamma) \\
    J^*(\pi_0) & \coloneqq \inf_{\gamma} J(\pi_0, \gamma).
\end{align}
We will also consider policies which obtain the above infima within some arbitrary threshold \(\epsilon > 0\), which we call \emph{near-optimal}, as follows:
\begin{definition}\label{definition:near-optimal}
    We say that a set of policies $\{\gamma\}$ depending on some parameter
set is  \emph{near-optimal} for the discounted distortion problem (respectively, average distortion problem) if  for any
$\epsilon > 0$, there is some choice of parameters such that the resulting policy  $\gamma$ satisfies \(J_\beta(\pi_0, \gamma) \le J_\beta^*(\pi_0) + \epsilon\)
(respectively, $J(\pi_0, \gamma) \le J^*(\pi_0) + \epsilon$).
\end{definition}

Note that in the zero-delay coding problem, we are usually concerned with the average distortion problem. However, it can be shown that as \(\beta \to 1\), a policy that is near-optimal for the discounted distortion problem is also near-optimal for the average distortion problem (see \cite[Theorem 7.3.6]{yuksel2023control}). Furthermore, the discounted distortion problem is generally easier to study from a reinforcement learning standpoint. Thus, we will target the discounted distortion problem throughout the majority of the paper and then make connections with the average distortion problem by taking \(\beta \to 1\).

For fixed \((x_{[0,t-1]}, m_{[0,t-1]}, m'_{[0,t-1]})\), consider the function \(\gamma(\cdot, x_{[0,t-1]}, m_{[0,t-1]}, m'_{[0,t-1]}) : \mathcal{X} \to \mathcal{M}\). Such a function (that is, a mapping from \(\mathcal{X}\) to \(\mathcal{M}\)) is called a \emph{quantizer}. We denote the set of all quantizers by \(\mathcal{Q}\). Thus we can view a policy \(\gamma\) as selecting a quantizer \(Q_t \in \mathcal{Q}\) based on the information \((X_{[0,t-1]}, M_{[0,t-1]}, M'_{[0,t-1]})\), then generating the channel input \(M_t\) as \(Q_t(X_t)\), as in~\cite{YukLinZeroDelay}.

Recall that we used \(O(m' | m)\) to denote our channel transition kernel. Let \(O_{Q}(m'| x)\) denote the kernel induced by a quantizer \(Q \in \mathcal{Q}\); that is, \(O_{Q}(m' | x) = O(m' | Q(x))\). Now let \(\psi \in \mathcal{P}(\mathcal{M'})\) be such that \(O_Q(\cdot | x) \ll \psi\) for all \(x \in \mathcal{X}, Q \in \mathcal{Q}\), where we use ``\(\ll\)'' to denote absolute continuity (that is, \(\psi(B) = 0 \implies O_Q(B | x) = 0\) for any Borel \(B \subset \mathcal{M'}\)). Since \(\mathcal{M'}\) is finite in our setup, we will take \(\psi\) to be the uniform measure on \(\mathcal{M'}\), but note that such measures also exists in uncountable setups for most practical channels. Then let \(g_Q(x,m') \coloneqq \frac{dO_Q}{d\psi}(x,m')\) be the Radon-Nikodym derivative of \(O_Q\) with respect to \(\psi\). In particular for a uniform \(\psi\), we have \(g_Q(x,m') = \left| \mathcal{M'} \right| O_Q(m' | x)\).


Also, let \(\pi_t, \overline{\pi}_t \in \mathcal{P}(\mathcal{X})\) be defined as
\begin{align}\label{predFiltDefn}
    \pi_t(\cdot)            & = P^\gamma_{\pi_0}(X_t \in \cdot | M'_{[0,t-1]}) \\
    \overline{\pi}_t(\cdot ) & = P^\gamma_{\pi_0}(X_t \in \cdot | M'_{[0,t]}),
\end{align}
recalling that \(X_0 \sim \pi_0\). We have dropped the \(\gamma\) for notational simplicity, but it should be noted that \(\pi_t\) and \(\overline{\pi_t}\) are policy-dependent. In the literature, \(\pi_t\) is called the \emph{predictor}, as it is the predictive probability of the next symbol $X_t$. \(\overline{\pi}_t\) is called the \emph{filter} due to its role in estimation and nonlinear filtering of partially observable systems \cite{van2007stochastic}. With a slight abuse of notation, we also let the source transition kernel \(T\) act as an operator on probability measures as follows:
\begin{align}
    T : \mathcal{P}(\mathcal{X}) & \to \mathcal{P}(\mathcal{X})             \\
    \pi(x)                       & \mapsto \sum_\mathcal{X}T(x' | x)\pi(x).
\end{align}

Then given \(\pi_0\), the above measures can be computed in a recursive manner as follows (see \cite[Proposition 3.2.5]{CappeHidden}).
\begin{align}
    \overline{\pi}_t(x) & = \frac{g_{Q_t}(x, M'_t)\pi_t(x)}{\sum_{\mathcal{X}} g_{Q_t}(x, M'_t) \pi_t(x)} \nonumber, \\
    \pi_{t+1}           & = T(\overline{\pi}_t) \label{eq:update} .
\end{align}

Using the above update equations, one can compute \(\pi_t\) given \((M'_{[0,t-1]}, Q_{[0,t-1]})\), so that policies of the form \(Q_t = \gamma_t(\pi_t)\) are valid. These policies form a special class.
\begin{definition}\cite{wood2016optimal}
    We say a policy \(\gamma = \{\gamma_t\}_{t\ge0}\) is of the \emph{Walrand-Varaiya type} if, at time \(t\), \(\gamma\) selects a quantizer \(Q_t = \gamma_t(\pi_t)\) and \(M_t\) is generated as \(M_t = Q_t(X_t)\). Such a policy is called \emph{stationary} if it does not depend on \(t\) (that is, \(\gamma_t = \overline{\gamma}\) for some \(\overline{\gamma}\) and all \(t \ge 0\)). The set of all stationary Walrand-Varaiya policies is denoted by \(\Gamma_{\text{WS}}\).
\end{definition}

The following are key results, originally from Walrand and Varaiya~\cite{WalrandVaraiya} for a finite time horizon and extended to the infinite-horizon case in~\cite{wood2016optimal}.

\begin{proposition}\cite[Proposition 2]{wood2016optimal}\label{proposition:discounted}
    For any \(\beta \in (0,1)\), there exists \(\gamma^* \in \Gamma_{\text{WS}}\) that solves the discounted distortion problem (that is, it minimizes~\eqref{eq:dis_cost}) for all priors \(\pi_0 \in \mathcal{P}(\mathcal{X})\).
\end{proposition}

\begin{proposition}\cite[Theorem 3]{wood2016optimal}
    There exists \(\gamma^* \in \Gamma_{\text{WS}}\) that solves the average distortion problem (that is, it minimizes~\eqref{eq:avg_cost}) for all priors \(\pi_0 \in \mathcal{P}(\mathcal{X})\).
\end{proposition}

\subsection{Regularity Properties of the Markov Decision Process}\label{section:MDPs}
Utilized in the above results is the fact that, under any \(\gamma \in \Gamma_{\text{WS}}\), the zero-delay coding problem can be viewed as a Markov decision process (MDP), which we now formally define.

\begin{definition}\label{definition:MDP}
    We define a \emph{Markov decision process} (MDP) as a 4-tuple \((\mathcal{Z},\mathcal{U},P,c)\), where:
    \begin{enumerate}
        \item \(\mathcal{Z}\) is the \emph{state space}, which we assume is Polish (a Borel subset of a complete, separable metric space).
        \item \(\mathcal{U}\) is the \emph{action space}, also Polish.
        \item \(P : \mathcal{Z} \times \mathcal{U} \to \mathcal{P}(\mathcal{Z})\) is the \emph{transition kernel}, such that \((z,u) \mapsto P(dz' | z,u)\).
        \item \(c : \mathcal{Z} \times \mathcal{U} \to [0,\infty)\) is the \emph{cost function}.
    \end{enumerate}
\end{definition}

The objective is to minimize \(J_\beta(z_0, \gamma) \coloneqq \mathbf{E}^{\gamma}_{z_0} \left[ \sum_{t=0}^\infty \beta^t c(Z_t, U_t)\right]\) (which we call the discounted cost problem) or \(J(z_0, \gamma) \coloneqq \limsup_{T \to \infty} \mathbf{E}^{\gamma}_{z_0} \left[ \frac{1}{T} \sum_{t=0}^{T-1} c(Z_t, U_t)\right]\) (which we call the average cost problem), over all \(\gamma\), where \(\gamma = (\gamma_t)_{t \ge 0}\) and \(U_t = \gamma_t(Z_{[0,t]}, U_{[0,t-1]})\).

\begin{proposition}\label{proposition:is_MDP}
    Under any \(\gamma \in \Gamma_{\text{WS}}\), the zero-delay coding problem is an MDP, where:
    \begin{enumerate}
        \item \(\mathcal{Z} = \mathcal{P}(\mathcal{X})\).
        \item \(\mathcal{U} = \mathcal{Q}\).
        \item \(P = P(d\pi'|\pi,Q)\) induced by the update equations in~\eqref{eq:update}.
        \item \(c(\pi,Q) = \sum_{\mathcal{M'}} \min_{\hat{x} \in \hat{\mathcal{X}}} \sum_{\mathcal{X}} d(x,\hat{x})O_Q(m'|x)\pi(x) \).
    \end{enumerate}
\end{proposition}

This follows directly from the update equations in~\eqref{eq:update} and the fact that, under any \(\gamma \in \Gamma_{\text{WS}}\), \(\pi_t\) completely determines \(Q_t\). The choice of \(c\) is due to the following result.

\begin{lemma}\label{lemma:cost}
    If an optimal decoder is used, the expected distortion at the encoder (that is, before sending \(M_t\)) is given by
    \begin{equation}\label{eq:cost}
        c(\pi_t,Q_t) = \sum_{\mathcal{M'}} \min_{\hat{x} \in \hat{\mathcal{X}}} \sum_{\mathcal{X}} d(x,\hat{x})O_{Q_t}(m'|x)\pi_t(x).
    \end{equation}
\end{lemma}

\begin{IEEEproof}
    Note that, under any $\gamma \in \Gamma_{\text{WS}}$ and given $\pi_0$, the optimal decoder can compute the filter $\overline{\pi}_t$, and thus the optimal decoder chooses \(\hat{X}_t\) according to
    \begin{equation}
        \argmin_{\hat{x}}\mathbf{E}\Bigl[ d(X_t, \hat{x}) | M'_{[0,t]} \Bigr] = \argmin_{\hat{x}} \sum_{x} d(x, \hat{x}) \overline{\pi}_t(x).
    \end{equation}
    By the update equations in~\eqref{eq:update}, we have
    \begin{equation}
        \overline{\pi}_t(x) = \frac{g_{Q_t}(x, M'_t)\pi_t(x)}{N_t(M'_t, Q_t)},
    \end{equation}
    where \(N_t(m', Q) = \sum_{\mathcal{X}} g_{Q}(x, m') \pi_t(x)\). Thus at the decoder, the expected distortion is given by
    \begin{equation}
        \min_{\hat{x}} \sum_{x} d(x, \hat{x}) \frac{g_{Q_t}(x, M'_t)\pi_t(x)}{N_t(M'_t, Q_t)}.
    \end{equation}

    However, at the encoder we must take the further expectation over \(M'_t\) (conditioned on \(M'_{[0,t-1]}\)), since we do not yet have access to \(M'_t\). Thus, at the encoder the expected distortion is
    \begin{align}
          & \sum_{m'} \min_{\hat{x}} \sum_{x} d(x, \hat{x}) \frac{g_{Q_t}(x, m')\pi_t(x)}{N_t(m', Q_t)} P^\gamma_{\pi_0}(M'_t = m' | M'_{[0,t-1]}) \\
        = & \sum_{m'} \min_{\hat{x}} \sum_{x} d(x, \hat{x}) \frac{g_{Q_t}(x, m')\pi_t(x)}{N_t(m', Q_t)} \sum_{x} P^\gamma_{\pi_0}(M'_t = m', X_t = x | M'_{[0,t-1]}) \\
        = & \sum_{m'} \min_{\hat{x}} \sum_{x} d(x, \hat{x}) \frac{g_{Q_t}(x, m')\pi_t(x)}{N_t(m', Q_t)} \sum_{x} P^\gamma_{\pi_0}(X_t = x | M'_{[0,t-1]}) P^\gamma_{\pi_0}(M'_t = m' | X_t = x, M'_{[0,t-1]}) \\
        = & \sum_{m'} \min_{\hat{x}} \sum_{x} d(x, \hat{x}) \frac{g_{Q_t}(x, m')\pi_t(x)}{N_t(m', Q_t)} \sum_{x} \pi_t(x) O_{Q_t}(m' | x)   \\
        = & \sum_{m'} \min_{\hat{x}} \sum_{x} d(x, \hat{x}) O_{Q_t}(m' | x)\pi_t(x),
    \end{align}
    where third equality holds by the definition of $\pi_t$ and since $P^\gamma_{\pi_0}(M'_t = m' | X_t = x, M'_{[0,t-1]}) = O_{Q_t}(m'|x)$ by the memoryless property of the channel. The last equality follows since $\sum_{x} \pi_t(x) O_{Q_t}(m' | x) = |\mathcal{M}'|^{-1} N_t(m',Q_t) $ and $g_{Q_t}(x, m') = |\mathcal{M}'| O_{Q_t}(m'|x)$.
\end{IEEEproof}

By this lemma, we have that the expected distortion at the encoder (assuming an optimal decoder), satisfies
\begin{equation}
    \mathbf{E}_{\pi_0}^{\gamma}\left[\sum_{t=0}^{T-1}c(\pi_t,Q_t)\right] = \mathbf{E}_{\pi_0}^{\gamma}\left[\sum_{t=0}^{T-1}d(X_t,\hat{X}_t)\right].
\end{equation}
Thus, this choice of \(c\) ensures that solving the MDP defined in Proposition~\ref{proposition:is_MDP} over all \(\gamma \in \Gamma_{\text{WS}}\) (that is, minimizing \(J_\beta(\pi_0, \gamma)\) or \(J(\pi_0, \gamma)\)) is equivalent to solving the zero-delay coding problem. Accordingly, we hereafter consider the discounted and average cost problems for this MDP (rather than the original discounted and average distortion problems). This allows us to use strategies from the literature of stochastic control; however, several complexities have been introduced:
\begin{itemize}
    \item While the source alphabet \(\mathcal{X}\) is finite, the state space of the MDP, \(\mathcal{P}(\mathcal{X})\), is uncountable. Furthermore, while our source process \((X_t)_{t \ge 0}\) is finite and irreducible (and hence has a unique invariant measure), there is no a priori reason for the MDP state process \((\pi_t)_{t \ge 0}\) to inherit these properties; in particular, irreducibility is too demanding.
    \item While we assume knowledge of the source transition kernel \(T\), the calculation of the transition kernel \(P(d\pi' | \pi, Q)\) is computationally demanding. 
\end{itemize}

Thus even if one can approximate the MDP state space \(\mathcal{P}(\mathcal{X})\) by some finite one, implementation of traditional MDP methods such as dynamic programming is difficult for this problem. This motivates the use of a reinforcement learning approach in which the calculation of these transition probabilities is unnecessary. We will cover this method in detail in Section~\ref{section:Q-learning}. Finally, although explicit computation of \(P(d\pi' | \pi, Q)\) is difficult, the following key structural result was obtained in~\cite{YukLinZeroDelay}.

\begin{lemma}\cite[Lemma 11]{YukLinZeroDelay}\label{lemma:weak_cts}
    The transition kernel \(P(d\pi' | \pi, Q)\) is weakly continuous. That is,
    \begin{equation}
        \int_{\mathcal{P}(\mathcal{X})} f(\pi')P(d\pi'|\pi,Q)
    \end{equation}
    is continuous in \((\pi,Q)\) for any continuous and bounded \(f : \mathcal{P}(\mathcal{X}) \to \mathbb{R}\).
\end{lemma}
Here, we endow \(\mathcal{P}(\mathcal{X})\) with the weak convergence topology and \(\mathcal{Q}\) with the Young topology (see \cite{YukLinZeroDelay}). Alternatively, since \(\mathcal{Q}\) is finite here the discrete topology would also suffice. We note that MDPs with weakly continuous transition kernels as above are often called \emph{weak Feller}.

\subsection{Filter and Predictor Stability}\label{section:pred_stability}
A key property that we use is \emph{filter/predictor} stability (recall from Definition \ref{predFiltDefn} that the predictor is given by \(\pi_t\) and the filter by \(\overline{\pi}_t\)).

\begin{definition}
    The total variation distance between two probability measures \(\mu, \nu\) defined over \(\mathcal{X}\) is given by
    \begin{equation}
        ||\mu - \nu||_{\text{TV}} \coloneqq \sup_{||f||_\infty \le 1} \left| \int_{\mathcal{X}} f(x)\mu(dx) - \int_{\mathcal{X}}f(x)\nu(dx) \right|,
    \end{equation}
    where the supremum is over all measurable real functions such that \(||f||_\infty = \sup_{x \in \mathcal{X}}|f(x)| \le 1\).
\end{definition}

Note that the total variation distance is equivalent to the \(L_1\) metric when \(\mathcal{X}\) is finite. Recall the update equations in~\eqref{eq:update} and note that they are sensitive to the value of \(\pi_0\). Accordingly, we use \(\pi_t^\mu\) to denote the predictor when \(\pi_0 = \mu\).

\begin{definition}\label{definition:stable}
    The predictor process \((\pi_t)_{t \ge 0}\) is \emph{stable in total variation almost surely} under a policy \(\gamma \in \Gamma_{\text{WS}}\) if, for all \(\mu, \nu \in \mathcal{P}(\mathcal{X})\) such that \(\mu \ll \nu\), we have
    \begin{equation}
        \lim_{t \to \infty} ||\pi_t^\mu - \pi_t^\nu||_{\text{TV}} = 0 \; P_\mu^\gamma\text{-a.s.},
    \end{equation}
    where $\pi^\mu_t$ (respectively, $\pi^\nu_t$) denotes the predictor with prior $\pi_0 = \mu$ (respectively, $\pi_0 = \nu$).
\end{definition}
That is, the predictor process is insensitive to its initialization. In some cases, we are interested in a different form of stability, which we define next.

\begin{definition}\label{definition:exp_stable}
    The predictor process \((\pi_t)_{t \ge 0}\) is \emph{exponentially stable in total variation in expectation} under a policy \(\gamma \in \Gamma_{\text{WS}}\) if there exists a coefficient \(\alpha \in (0,1)\) such that for all \(\mu, \nu \in \mathcal{P}(\mathcal{X})\) (with \(\mu \ll \nu\)) and for all \(t \ge 0\), we have
    \begin{equation}
        \mathbf{E}_\mu^\gamma \left[ ||\pi_{t+1}^\mu - \pi_{t+1}^\nu||_{\text{TV}} \right] \le \alpha \mathbf{E}_\mu^\gamma \left[ ||\pi_t^\mu - \pi_t^\nu||_{\text{TV}} \right].
    \end{equation}
\end{definition}

We can make equivalent definitions for the filter process; in fact, the problem of filter stability (in various senses) is a classical  problem in probability and statistics, where it is typically established in two ways: (i) The transition kernel of the underlying state is in some sense \emph{sufficiently ergodic}, so that regardless of the observations, the filter process inherits this ergodicity and forgets its prior over time. (ii) The observations are in some sense \emph{sufficiently informative}, so that, regardless of the prior, the filter process tracks the true state process. For a detailed review of these filter stability methods, see~\cite{Chigansky}. However, we will need slightly more general results in our case, since it is usually assumed in the filter stability problem that the observation kernel is time-invariant; here, \(O_{Q_t}\) depends on \(Q_t\) and hence changes with time, and accordingly additional analysis is needed.

\subsection{A Note on MDP Notation}\label{subsection:MDP_notation}
The notation in the following sections can get intricate, so we first introduce some additional notation to be used in the context of MDPs. Some of the concise discussion below will be expanded upon and made more specific in the following sections.
\begin{enumerate}
    \item When discussing \emph{approximations} of an MDP state we use a caret symbol. For example, we use \(\hat{Z}_t\) to denote an approximation of \(Z_t\). Accordingly we use \(\hat{\gamma}\) to denote a policy that maps \(\hat{Z}_t\) to \(U_t\) and we use \(\hat{J}_\beta(\hat{z}_0, \hat{\gamma})\) to be some appropriately defined discounted cost under that policy. We have to be careful with these definitions as the approximate sequence \((\hat{Z}_t)_{t \ge 0}\) may not form a Markov chain, but we cover those technicalities later. Furthermore, when the approximation depends on some parameter \(N\), we denote such a policy by \(\hat{\gamma}_N\); when used in this way, we assume the policy is stationary, and thus the subscript \(N\) should not be confused with a time index. In Section~\ref{section:finite window}, $N$ determines the window length, and in Section~\ref{subsection:NN approximation}, it determines the ``resolution'' of quantization.
    \item To denote \emph{optimality} we use an asterisk symbol. For example, we denote \(J_\beta^*(z_0) \coloneqq \min_{\gamma} J_\beta(z_0, \gamma)\), and we denote a policy achieving this minimum by \(\gamma^*\).
    \item Finally, when we \emph{extend} an approximation back to the original state space, we use a tilde symbol. For example, if we wished to take \(\hat{\gamma}\) and appropriately modify it to take \(Z_t\) as an input rather than \(\hat{Z}_t\), we would call the resulting policy \(\widetilde{\gamma}\).
\end{enumerate}
Note that this notation is presented in order of operation. For example, \(\hat{J}_\beta^*(\hat{z}_0)\) should be taken as the following: first we take an approximation of the MDP, then find the optimal discounted cost of this approximation. It is not the approximation of the optimal discounted cost \(J_\beta^*(z_0)\), although under some conditions to be presented later it may be interpreted as such. Accordingly, \(\widetilde{J}_\beta^*(z_0)\) should be taken as the extension of \(\hat{J}_\beta^*(\hat{z}_0)\) to all of \(\mathcal{Z}\), not as the optimal discounted cost for some ``extended'' MDP.

\section{Sliding Finite Window Approximation of the Belief MDP}\label{section:finite window}
Here we describe how to approximate \(\pi_t\) using a sliding finite window. The analysis in this section is inspired by~\cite{kara2021convergence}, which used a similar construction to study sliding finite window policies for partially observed Markov decision processes (POMDPs). We note that although the proof methods are similar, our setup differs from \cite{kara2021convergence} in two ways: 1) our problem is not strictly a POMDP as studied in \cite{kara2021convergence}, because the induced observation kernel changes with the choice of quantizer $Q_t$, 2) accordingly, our belief state is the predictor rather than the filter, and thus the Dobrushin coefficient terms we arrive at are different with a tailored analysis.

First, we must define a slightly different MDP than that defined in Proposition~\ref{proposition:is_MDP}. Fix some window size \(N \in \mathbb{Z}_+\). Recall the channel outputs \((m'_t)_{t \ge 0}\) and the quantizers \((Q_t)_{t \ge 0}\). We define the following:
\begin{align}
    I_t & \coloneqq (M'_{[t-N, t-1]}, Q_{[t-N, t-1]}) \\
    W_t & \coloneqq
    (\pi_{t-N}, I_t).
\end{align}

Note that we can compute \(\pi_t\) given \(w_t\) by applying the update equations in~\eqref{eq:update} \(N\) times. Denote this mapping by
\begin{align}
    \varphi : \mathcal{W} & \to \mathcal{P}(\mathcal{X}) \\
    W_t                   & \mapsto \pi_t
\end{align}
where \(\mathcal{W} = \mathcal{P}(\mathcal{X}) \times (\mathcal{M'})^N \times \mathcal{Q}^N\), endowed with the product topology, where we use the weak convergence topology on \(\mathcal{P}(\mathcal{X})\) and standard coordinate topologies on \(\mathcal{M'}\) and \(\mathcal{Q}\).

We call \(W_t\) the sliding finite window belief term, and call policies of the form \(Q_t = \gamma_t(W_t)\) \emph{sliding finite window belief} policies (with window size \(N\)). If it does not depend on \(t\), we call it stationary. Denote the set of all stationary sliding finite window belief policies by \(\Gamma_{FS}\).

\emph{Remark:} This approach assumes that we start at time \(t \ge N\); although for discounted MDPs the first \(N\) steps may be significant, for the zero-delay coding problem we are interested in taking \(\beta \to 1\), so these first \(N\) steps will not be crucial. Accordingly, we assume that \(N\) steps  have already been completed with some arbitrary \(\gamma \in \Gamma_{\text{WS}}\). For notational simplicity, we assume that these steps have occurred from \(t = -N, \ldots, -1\) and thus the process starts from \(w_0\) (and the prior would now be \(\pi_{-N}\)).

\subsection{The Sliding Finite Window Belief MDP}
This sliding finite window belief construction inherits the MDP properties of the original setup. Indeed, it is straightforward to show that the process \((W_t)_{t \ge 0}\) is controlled Markov, with control \((Q_t)_{t \ge 0}\). That is, for all \(t \ge 0\),
\begin{equation}
    P(W_{t+1} \in \cdot | W_{[0,t]}, Q_{[0,t]}) = P(W_{t+1} \in \cdot | W_t, Q_t).
\end{equation}

\begin{proposition}\label{proposition:is_MDP_2}
    Under any \(\gamma \in \Gamma_{FS}\), the zero-delay coding problem is an MDP, where:
    \begin{enumerate}
        \item \(\mathcal{Z} = \mathcal{W}\).
        \item \(\mathcal{U} = \mathcal{Q}\).
        \item \(P = P(dw'|w,Q)\).
        \item \(c(w,Q) = \sum_{\mathcal{M'}} \min_{\hat{x} \in \hat{\mathcal{X}}} \sum_{\mathcal{X}} d(x,\hat{x})O_Q(m'|x)\varphi(w)(x) \).
    \end{enumerate}
\end{proposition}
Note that the cost function is exactly the cost function we had in Proposition~\ref{proposition:is_MDP}, by simply replacing \( \pi = \varphi(w)\). Then an analog to Lemma~\ref{lemma:cost} holds, and we indeed have that solving the MDP from Proposition~\ref{proposition:is_MDP_2} is equivalent to solving the zero-delay coding problem. That is, we can equivalently consider \(J^*_\beta(w_0) = \inf_{\gamma \in \Gamma_{FS}} J_\beta(w_0, \gamma)\). The next proposition follows immediately from Proposition~\ref{proposition:discounted} and the fact that \(\pi_t = \varphi(W_t)\).
\begin{proposition}
    For any \(\beta \in (0,1)\) and \(N \in \mathbb{Z}_+\), there exists \(\gamma^* \in \Gamma_{FS}\) that solves the discounted cost problem; that is, one that satisfies, for all \(w_0 \in \mathcal{W}\),
    \[J_\beta(w_0,\gamma^*) = J_\beta^*(w_0).\]
\end{proposition}

\subsection{Sliding Finite Window Approximation}
The above representation is not particularly useful, as it still requires one to compute \(\pi_{t-N}\). Instead, fix the first coordinate to \(\zeta\) (the invariant distribution of the source) and let
\begin{align}\label{eq:pi_hat}
    \hat{W}_t   & = (\zeta, I_t)         \\
    \hat{\pi}_t & = \varphi(\hat{W}_t).
\end{align}
That is, we obtain \(\hat{\pi}_t\) by applying the update equations \(N\) times, but starting from an incorrect (fixed) prior \(\zeta\). Equivalently,
\begin{align}
    \pi_t(x)       & = P^\gamma_{\pi_{t-N}} (X_t = x | M'_{[t-N,t-1]}, Q_{[t-N,t-1]}) \\
    \hat{\pi}_t(x) & = P^\gamma_{\zeta} (X_t = x | M'_{[t-N,t-1]}, Q_{[t-N,t-1]}).
\end{align}
The key idea, which will be discussed in detail later, is that under predictor stability the correct predictor \(\pi_t = \varphi(w_t)\) and the incorrect predictor \(\hat{\pi}_t = \varphi(\hat{w}_t)\) will be close for large enough \(N\), since the predictor will be insensitive to the prior.

The benefits of such an approximation are evident: rather than deal with all of \(\mathcal{W}\), which is uncountable due to \(\mathcal{P}(\mathcal{X})\), we only have to deal with the finite set \(\mathcal{W}_N \coloneqq \{\zeta\} \times (\mathcal{M'})^N \times \mathcal{Q}^N\). Furthermore, we no longer need to compute \(\pi_{t-N}\), which can save significant computation resources especially when the relevant alphabets are large.

Consider the following transition kernel,
\begin{equation}
    P_N(\hat{w}_1 | \hat{w}, Q) \coloneqq P(\mathcal{P}(\mathcal{X}), i_1 | \hat{w}, Q),\label{eq:eta_hat}
\end{equation}
where \(P\) is the transition kernel of the sliding finite window belief MDP and \(\hat{w}_1 = (\zeta, i_1)\), and cost function
\begin{equation}
    c_N(\hat{w}, Q) \coloneqq \sum_{\mathcal{M'}} \min_{\hat{x} \in \hat{\mathcal{X}}} \sum_{\mathcal{X}} d(x,\hat{x})O_Q(m'|x)\hat{\pi}(x).\label{eq:c_hat}
\end{equation}

Then, our approximate MDP becomes \(\text{MDP}_N = (\mathcal{W}_N, \mathcal{Q}, P_N, c_N)\). As in previous sections, denote the discounted cost for this MDP under a given policy $\hat{\gamma}_N$ (from \(\mathcal{W}_N\) to \(\mathcal{Q}\)) by \(\hat{J}_\beta(\hat{w}_0, \hat{\gamma}_N)\), and the optimal discounted cost by \(\hat{J}^*_\beta(\hat{w}_0)\), with minimizing policy \(\hat{\gamma}_N^*\). The relevant extensions (which are simply obtained by making the previous functions constant over \(\mathcal{P}(\mathcal{X})\)) are then \(\widetilde{J}_\beta^*(w_0)\) and \(\widetilde{\gamma}^*_N\). The following is a key loss term:
\begin{align}\label{eq:loss}
    L^N_t \coloneqq \sup_{\gamma \in \Gamma_{\text{WS}}} \mathbf{E}^\gamma_{\pi_{t-N}} \left[ ||\pi_t - \hat{\pi}_t||_{\text{TV}} \right].
\end{align}
We now present our main results for this approximation scheme, which give a bound on the performance loss when using the given window length \(N\). Note that here we take an expectation with respect to some policy acting on the previous \(N\) steps (and hence generates \(W_0\)), and with respect to some prior \(\pi_{-N}\). Also, define $||d||_\infty \coloneqq \max_{x, \hat{x}}d(x, \hat{x})$, where we recall that $d$ is the distortion measure for the zero-delay coding problem.

\begin{theorem}\label{theorem:value_convergence} For any \(\gamma \in \Gamma_{\text{WS}}\) acting on \(N\) time steps to generate \(W_0\) and any prior \(\pi_{-N} \in \mathcal{P}(\mathcal{X})\), we have
    \begin{align}
        \mathbf{E}^\gamma_{\pi_{-N}} \left[ \left| \widetilde{J}_\beta^*(W_0) - J^*_\beta(W_0) \right| \right] \le \frac{||d||_\infty}{1-\beta}\sum_{t=0}^\infty \beta^t L^N_t.
    \end{align}
\end{theorem}
That is, we can bound the difference in the ``extended'' approximate optimal cost and the true optimal cost by a predictor-distance term.
\begin{theorem}\label{theorem:finite window_near_optimal}
    For any \(\gamma \in \Gamma_{\text{WS}}\) acting on \(N\) time steps to generate \(w_0\) and any prior \(\pi_{-N} \in \mathcal{P}(\mathcal{X})\), we have
    \begin{align}
        \mathbf{E}^\gamma_{\pi_{-N}} \left[ \left| J_\beta(W_0, \widetilde{\gamma}^*_N) - J^*_\beta(W_0) \right| \right] \le \frac{2||d||_\infty}{1-\beta}\sum_{t=0}^\infty \beta^t L^N_t.
    \end{align}
\end{theorem}
That is, we can bound the sub-optimality of the ``extended'' approximate optimal policy $\widetilde{\gamma}^*_N$ by a predictor-distance term. Thus, if we can show that this distance goes to $0$ as $N$ grows, and in particular if we can determine a certain rate of convergence, we can get a similar rate of convergence to optimality for our policy $\widetilde{\gamma}^*_N$.
The proofs for Theorems~\ref{theorem:value_convergence} and~\ref{theorem:finite window_near_optimal} are given in the Appendix.

\subsection{Bounds on the Loss Term}
The loss term \(L^N_t\) in the previous theorems is related to the question of predictor stability (recall Section \ref{section:pred_stability}). Indeed, the term within the supremum is exactly
\begin{align}
    \mathbf{E}^\gamma_\mu \left[ ||\pi^\mu_t - \pi^\nu_t||_{\text{TV}} \right]\label{eq:predictor-stability}
\end{align}
when \(\mu = \pi_{t-N}\) and \(\nu = \zeta\), and under some \(\gamma \in \Gamma_{\text{WS}}\). Thus bounding this term over all \(\gamma\) will give us a bound on \(L^N_t\). We note that any notion of predictor stability could be used to give a bound on $L^N_t$ (and thus on the performance of $\widetilde{\gamma}^*_N$). In the following section, we give one such condition which is sufficient (but by no means necessary) to apply the above theorems.

\subsection*{Dobrushin Coefficient Conditions}
The following results are inspired by the analysis in~\cite{mcdonald2020exponential}, which uses joint contraction properties of the state and observation kernels to bound~\eqref{eq:predictor-stability}. First we introduce some notation. For standard Borel spaces \(\mathcal{A}_1, \mathcal{A}_2\) and some kernel \(K : \mathcal{A}_1 \to \mathcal{P}(\mathcal{A}_2)\), we define the Dobrushin coefficient as
\begin{align}
    \delta(K) \coloneqq \inf \sum_{i=1}^n \min(K(B_i | x), K(B_i | y)),
\end{align}
where the infimum is over \(x, y \in \mathcal{A}_1\) and all partitions \(\{B_i\}_{i=1}^n\) of \(\mathcal{A}_2\). In particular, for finite spaces, the Dobrushin coefficient is equivalent to summing the minimum elements between every pair of rows, then taking the minimum of these sums. For example, take
\begin{equation}
    K = \begin{pmatrix}
        \frac{1}{2} & \frac{1}{2} & 0           & 0           \\
        \frac{1}{3} & \frac{1}{3} & \frac{1}{3} & 0           \\
        \frac{1}{3} & \frac{1}{3} & 0           & \frac{1}{3} \\
        \frac{1}{4} & \frac{1}{4} & \frac{1}{4} & \frac{1}{4}
    \end{pmatrix}.
\end{equation}
Between the first and second rows, the sum of the minimum elements gives~\(\frac{2}{3}\), between the third and fourth gives~\(\frac{3}{4}\), etc. One can verify that the minimum of such sums is \(\frac{1}{2}\), so \(\delta(K) = \frac{1}{2}\) (note that \(\delta(K) \le 1\) by definition).

The following is then a counterpart of~\cite[Theorem 3.6]{mcdonald2020exponential}. Note that in our case, the channel is not time-invariant, unlike in~\cite{mcdonald2020exponential}, but the analysis follows similarly. The proof is provided in the Appendix.

\begin{theorem}\label{theorem:dobrushin}
    For any \(\gamma \in \Gamma_{\text{WS}}\) and for any \(\mu \ll \nu\), we have
    \begin{align}
         & \mathbf{E}_{\mu}^\gamma \left[ ||\pi^\mu_{t+1} - \pi^\nu_{t+1}||_{\text{TV}} \right] \le (1-\delta(T))(2 - \tilde{\delta}(O))\mathbf{E}_{\mu}^\gamma \left[ ||\pi^\mu_t - \pi^\nu_t||_{\text{TV}} \right],
    \end{align}
    where \(\tilde{\delta}(O) = \min_{Q \in \mathcal{Q}}(\delta(O_Q))\).
\end{theorem}

We can arrive at a simpler bound by using \(\delta(O)\) directly rather than \(\tilde{\delta}(O)\). To see this, note that for a given quantizer \(Q\), the kernel \(O_Q(m' | x) = O(m' | Q(x))\) only contains rows from the kernel \(O\), thus \(\delta(O) \le \delta(O_Q)\) for all \(Q\). Thus we arrive at the following corollary.

\begin{corollary}\label{corollary:dobrushin}
    Assume \(\alpha \coloneqq (1-\delta(T))(2 - \delta(O)) < 1\). Then for any \(\gamma \in \Gamma_{\text{WS}}\) and for any \(\mu \ll \nu\), we have
    \begin{align}
        \mathbf{E}_\mu^\gamma \left[ ||\pi^\mu_{t+1} - \pi^\nu_{t+1}||_{\text{TV}} \right] \le \alpha \mathbf{E}_{\mu}^\gamma \left[ ||\pi^\mu_t - \pi^\nu_t||_{\text{TV}} \right].
    \end{align}
    That is, the predictor process is exponentially stable in total variation in expectation. Furthermore, if \(\delta(T) > \frac{1}{2}\), then the above is true with \(\alpha = 1 - \delta(T)\) regardless of the channel \(O\).
\end{corollary}

Applying this to the \(L^N_t\) term, we have
\begin{align}
    L^N_t & = \sup_{\gamma \in \Gamma_{\text{WS}}} \mathbf{E}^\gamma_{\pi_{t-N}} \left[ ||\pi_t - \hat{\pi}_t||_{\text{TV}} \right] \nonumber \\
          & \le \alpha^N ||\pi_{t-N} - \pi'||_{\text{TV}} \le 2\alpha^N.\label{eq:loss_bound}
\end{align}

The requirement \((1-\delta(T))(2 - \delta(O)) < 1\) places conditions on the source dynamics and the channel. A universal condition to make the result applicable over all channels is obtained when one considers the special case where the channel is noiseless: in this case, we have \(\delta(O) = 0\) and \(\delta(T) > \frac{1}{2}\) is to hold. This is not surprising given the nature of Dobrushin-type conditions. At a high level, the Dobrushin coefficient tells us how similar the conditional measures \(O_Q(m' | x)\) and \(O_Q(m' | y)\) are for different \(x, y \in \mathcal{X}\). The more similar they are, the closer the coefficient is to \(1\).


Combining Theorem~\ref{theorem:finite window_near_optimal} and the bound in~\eqref{eq:loss_bound}, we obtain the following result.
\begin{corollary}\label{corollary:dobrushin_finite_window}
    Assume \(\alpha \coloneqq (1-\delta(T))(2 - \delta(O)) < 1\). Then for any \(\gamma \in \Gamma_{\text{WS}}\) which acts on \(N\) time steps to generate \(W_0\) and any prior \(\pi_{-N}\), we have
    \begin{align}
        \mathbf{E}^{\gamma}_{\pi_{-N}} \left[ \left| J_\beta(W_0, \widetilde{\gamma}^*_N) - J^*_\beta(W_0) \right| \right] \le \frac{4||d||_\infty}{(1-\beta)^2}\alpha^N.
    \end{align}
\end{corollary}

\medskip

\section{Q-learning: Convergence to Near-Optimality}\label{section:Q-learning}
We now turn our attention to actually calculating an optimal policy for the approximate MDP in Section~\ref{section:finite window}, which will be near-optimal for the zero-delay coding problem for sufficiently large parameter \(N\). That is, the remaining problem is to find for some large $N$ the $\hat{\gamma}^*_N$ which minimizes
\begin{align}
    \mathbf{E}_{\pi_{-N}}^{\hat{\gamma}^*_N}\left[\sum_{t \ge 0} \beta^t c_N(\hat{w}_t, Q_t)\right].
\end{align}
We can then extend this policy in an appropriate manner to obtain a near-optimal policy $\widetilde{\gamma}^*$ for
\begin{align}
    \mathbf{E}_{\pi_0}^{\widetilde{\gamma}^*}\left[ \sum_{t \ge 0} \beta^t c(\pi_t, Q_t)\right],
\end{align}
which was equivalent to the discounted distortion in \eqref{eq:dis_cost}.

A tempting strategy from stochastic control theory would be the \emph{value iteration} method (see \cite[Lemma 4.2.8]{Hernandez}), which iteratively computes
\begin{equation}
    V_t(\hat{w}) = \min_{\mathcal{Q}}\left( c_N(\hat{w},Q) + \beta \int_{\mathcal{W}_N} V_{t-1}(\hat{w}')P_N(d\hat{w}' | \hat{w}, Q) \right),
\end{equation}
for all \(\hat{w} \in \mathcal{W}_N\) and with \(V_0(\cdot) \equiv 0\). However, if we attempt to perform this method for the approximate zero-delay MDP, it becomes very analytically complicated. In particular, although we know the transition matrix \(T\) of the source \((X_t)_{t \ge 0}\), it is not clear how to compute \(P_N\) for the approximate MDP. This motivates the use of a reinforcement learning algorithm in which this transition kernel is not computed directly and the optimal policy is learned empirically. We study a version of Q-learning, originally from~\cite{Watkins} for finite sets and generalized for weak Feller MDPs (along with rigorous convergence and near-optimality guarantees) in~\cite{KSYContQLearning} and~\cite{kara2021convergence}. We first define the following update equations, where \(V_t\) and \(\alpha_t\) are both functions from \(\mathcal{Z} \times \mathcal{U}\) to \(\mathbb{R}_+\), and \(C_t \in \mathbb{R}_+\):
\begin{align}
    V_{t+1}(z,u)     & = V_t(z,u) \quad \text{for all } (z,u) \neq (Z_t,U_t), \nonumber                                                                                              \\
    V_{t+1}(Z_t,U_t) & = \left( 1 - \alpha_t(Z_t,U_t) \right)V_t(Z_t,U_t) + \alpha_t(Z_t,U_t)\left[ C_t + \beta \min_{v \in \mathcal{U}}V_t(Z_{t+1}, v) \right], \label{eq:V_update} \\
    \alpha_t(z,u)    & = \frac{1}{1 + \sum_{k=0}^t\mathbf{1}(Z_k = z,U_k = u)} \nonumber,
\end{align}
for some arbitrary \(V_0\). Note that in the literature, \(Q\) is often used instead of \(V\); here we reserve \(Q\) for quantizers. Our first algorithm is simply data collection - it samples \(X_t\), \(Q_t\), and \(M'_t\) to be used in our Q-learning iterations. Note that although the data collection and the Q-learning are presented as two separate algorithms for clarity, one could run the two concurrently: as we will see, all of the information needed to compute \(V_t\) is present at time \(t\) in the data collection.

\bigskip
\medskip

\noindent\underline{\textbf{Algorithm 1: Data Collection}}
\begin{algorithmic}[1]
    \Require initial distribution \(\pi_0\), transition kernel \(T\), quantizer set \(\mathcal{Q}\), channel kernel \(O\)
    \State Sample \(X_0 \sim \pi_0\)
    \State Choose \(Q_0\) uniformly from \(\mathcal{Q}\)
    \State \(M_0 = Q_0(X_0)\)
    \State Sample \(M'_0 \sim O(\cdot | M_0)\)
    \For{\(t \ge 1\)}
    \State Sample \(X_t \sim T(\cdot | X_{t-1})\)
    \State Choose \(Q_t\) uniformly from \(\mathcal{Q}\)
    \State \(M_t = Q_t(X_t)\)
    \State Sample \(M'_t \sim O(\cdot | M_t)\)
    \EndFor
\end{algorithmic}

This gives us a sequence \((Q_t, M'_t)_{t \ge 0}\). We then apply the iterations in~\eqref{eq:V_update} with appropriate choices for \(Z_t\), \(U_t\) and \(C_t\). With an abuse of notation, let us shift the time index of the sampled data by \(N\) to obtain \((Q_t, M'_t)_{t \ge -N}\). Then apply the sliding finite window scheme from Section~\ref{section:finite window} to obtain \((\hat{W}_t)_{t \ge 0}\). Then we define \(Z_t\), \(U_t\) and \(C_t\) as follows:
\begin{align}
    Z_t & \coloneqq \hat{W}_t            \\
    U_t & \coloneqq Q_t                  \\
    C_t & \coloneqq c_N(\hat{W}_t, Q_t).
\end{align}
\emph{Remark:} When this method is used in \cite{kara2021convergence} for POMDPs, the realizations of the ``true'' cost \(c(z, u)\) is used, since it is not assumed that we know the exact form of the cost function. However, in our case we wish to avoid computing \(\pi_t\) since it is not needed for the approximation, and computing \(c_N\) is straightforward via \eqref{eq:c_hat}. Accordingly, we use the ``approximate'' cost \(c_N(\hat{W}_t, Q_t)\).

\begin{theorem}\label{theorem:Qlearning2}
    The sequence \((V_t)_{t \ge 0}\) converges almost surely to a limit \(V^*\). Furthermore, consider the policy
    \begin{equation}
        \hat{\gamma}^*_N(\hat{w}) = \argmin_{Q \in \mathcal{Q}} V^*(\hat{w},Q),
    \end{equation}
    and let \(\widetilde{\gamma}^*_N\) be the extension of \(\hat{\gamma}^*_N\) to \(\mathcal{W}\) by making it constant over the belief coordinate, as in Theorem~\ref{theorem:finite window_near_optimal}. Then for any \(\gamma\) which acts on \(N\) time steps to generate \(W_0\) and any prior \(\pi_{-N}\),
    \begin{equation}
        \mathbf{E}^\gamma_{\pi_{-N}} \left[ J_\beta(W_0, \widetilde{\gamma}^*_N) - J^*_\beta(W_0) \right] \le \frac{4||d||_\infty}{(1-\beta)^2}\alpha^N,
    \end{equation}
    where $\alpha = (1-\delta(T))(2 - \delta(O))$.
\end{theorem}

The proof for Theorem \ref{theorem:Qlearning2} is given in the Appendix. Note that convergence relies on recent results in \cite{kara2023qlearning}; although convergence for finite MDPs is classical from \cite{Watkins}, when applied to approximate MDPs (as we do here), convergence is not trivial since the underlying true MDP dynamics differ from that of the approximate MDP. Accordingly we must verify certain ergodic behavior of the MDP to apply the results in \cite{kara2023qlearning}.

\section{An Alternative Approximation Scheme and Comparison}\label{section:comparison}
\subsection{Nearest Neighbor Approximation}\label{subsection:NN approximation}
Another approximation scheme for the zero-delay coding problem was proposed in \cite{creggZeroDelayNoiseless}, along with a Q-learning algorithm to compute the corresponding near-optimal policy. Although this scheme was only formally analyzed in the case of a noiseless channel, the results go through in the case of a noisy channel with some modifications. Accordingly, we present the noisy channel analogs of the results in \cite{creggZeroDelayNoiseless} without proof, but we do provide a detailed comparison of the two schemes.

First, we approximate $\pi_t$ using a nearest neighbor quantization scheme with a finite number of bins. This approximation can be done efficiently using~\cite[Algorithm 1]{Reznik}, which approximates \(\mathcal{P}(\mathcal{X})\) by the following finite set,
\begin{equation}
    \mathcal{P}_N(\mathcal{X}) = \left\{ \hat{\pi} \in \mathcal{P}(\mathcal{X}) : \hat{\pi} = \left[ \frac{k_1}{N}, \ldots, \frac{k_{|\mathcal{X}|}}{N} \right], k_i = 0, \ldots, N, i = 1, \ldots, |\mathcal{X}| \right\},\label{eq:P_N}
\end{equation}
where $\sum_{i=1}^N k_i=N$. Note that in the literature, such distributions are also called types or empirical distributions~\cite{Csiszar1998}. We let \(\hat{\pi}\) be the nearest neighbor (in Euclidean distance) of \(\pi\) in \(\mathcal{P}_N(\mathcal{X})\), and we clearly have that \(\max_{\pi}d(\pi, \hat{\pi}) \to 0\) as \(N \to \infty\).

Similarly to Section \ref{section:finite window}, one can then define a transition kernel and cost function over $\mathcal{P}_N(\mathcal{X})$. We omit the details, but they are essentially the averages of the true transition kernel $P(d\pi' | \pi, Q)$ and cost function $c(\pi,Q)$ over the nearest neighbor bins defined by \eqref{eq:P_N}, with respect to an appropriate reference measure. We denote the resulting approximate MDP by $\text{MDP}_N = (\mathcal{P}_N, \mathcal{Q}, P_N, c_N)$.

\begin{corollary}\cite[Theorem 1]{creggZeroDelayNoiseless}\label{corollary:quantization_near_optimal}
    Let \(\hat{\gamma}^*_N \in \Gamma_{\text{WS}}\) be optimal for \(\text{MDP}_N\) and let \(\widetilde{\gamma}^*_N(\pi) = \hat{\gamma}^*_N(\hat{\pi})\). Then for all \(\pi_0 \in \mathcal{P}(\mathcal{X})\) and \(\beta \in (0,1)\),
    \begin{equation}
        \lim_{N \to \infty} | J_\beta(\pi_0, \widetilde{\gamma}^*_N) - J_\beta^*(\pi_0) | = 0.
    \end{equation}
    That is, \(\widetilde{\gamma}^*_N\) is near-optimal for the zero-delay coding problem for sufficiently large \(N\).
\end{corollary}

The analogous Q-learning result is as follows, using the sampled data obtained via Algorithm 1. Let \((\pi_t)_{t \ge 0}\) be the sequence of predictors resulting from the sampled data \((Q_t, M'_t)_{t \ge 0}\) (with prior \(\pi_0\)). Recall that this can be computed recursively using the equations in~\eqref{eq:update}. Apply the nearest neighbor approximation scheme onto $\mathcal{P}_N(\mathcal{X})$ to obtain \((\hat{\pi}_t)_{t \ge 0} \subset \mathcal{P}_N(\mathcal{X})\). Then, recalling the iterations in \eqref{eq:V_update}, we define \(Z_t\), \(U_t\) and \(C_t\) as follows:
\begin{align}
    Z_t & \coloneqq \hat{\pi}_t                                          \\
    U_t & \coloneqq Q_t                                                  \\
    C_t & \coloneqq c(\pi_t, Q_t) \quad \text{(recall~\eqref{eq:cost})}.
\end{align}

\begin{theorem}\cite[Theorem 1]{creggZeroDelayNoiseless}\label{theorem:Qlearning1}
    Assume that \(\pi_0\) in Algorithm 1 is the unique invariant distribution of the source, \(\zeta\). Then \(V_t\) converges almost surely to a limit \(V^*\). Furthermore, consider the policy
    \begin{equation}
        \hat{\gamma}^*_N(\hat{\pi}) = \argmin_{Q \in \mathcal{Q}} V^*(\hat{\pi},Q),
    \end{equation}
    and let \(\widetilde{\gamma}^*_N\) be the extension of \(\hat{\gamma}^*_N\) to \(\mathcal{P}(\mathcal{X})\) by making it constant over the belief quantization bins, as in Corollary~\ref{corollary:quantization_near_optimal}. Then as \(N \to \infty\),
    \begin{equation}
        J_\beta(\zeta, \widetilde{\gamma}^*_N) \to J^*_\beta(\zeta).
    \end{equation}
    That is, \(\widetilde{\gamma}^*_N\) is near-optimal for the discounted distortion zero-delay coding problem, provided that the source starts from its invariant distribution.
\end{theorem}

\subsection{Reinforcement Learning Theoretic Comparison of the Two Schemes}
In this section, we provide a detailed comparison noting explicit benefits and drawbacks of each of the two learning schemes to attain near-optimality. In particular, we emphasize the easier implementation and lower computational complexity of our proposed sliding window method, although there are some cases where the belief quantization method of \cite{creggZeroDelayNoiseless} may be preferable.
\begin{itemize}
    \item \emph{Assumptions for convergence of Q-learning:} In Theorem~\ref{theorem:Qlearning1}, one can see that the convergence of \(V_t\) depends on the initial distribution used during learning (in particular, the \(\pi_0\) in Algorithm 1 must be the source's unique invariant distribution \(\zeta\)). Conversely, in Theorem~\ref{theorem:Qlearning2}, the convergence of \(V_t\) happens regardless of the initial distribution used during learning (that is, the \(\pi_0\) in Algorithm 1 can be arbitrary).
    
    \item \emph{Conditions on initialization during implementation:} For the belief quantization scheme, the policy from Q-learning is near-optimal when \(\pi_0 = \zeta\). We note that in the special case of the noiseless channel, valid initializations also include $\pi_0 = \delta_x$, i.e., deterministically starting from any $x \in \mathcal{X}$. This can be shown using \cite[Corollary 2]{creggZeroDelayNoiseless} and by noting that in this special case the Dirac masses $\delta_x$ are recurrent. However, in the noisy channel case, this argument fails and we require \(\pi_0 = \zeta\); conversely, for the sliding finite window scheme, the policy obtained in Theorem~\ref{theorem:Qlearning2} is near-optimal for almost every initial window.
   
    \item \emph{Filter stability assumptions:} Because the predictor stability required for the sliding finite window scheme is rather strong, the near-optimality of approximations was only established for those source-channel combinations satisfying the Dobrushin coefficient conditions. However, we note that the Dobrushin coefficent conditions are not necessary; any filter stability result which gives a bound on the loss term in Theorem \ref{theorem:finite window_near_optimal} could be used. To this end, Hilbert metric approaches may also be relevant, for example those in \cite{le2004stability, demirci2024refinedboundsnearoptimality}.
    
    \item \emph{Computational complexity and need for Bayesian updates:} In the belief quantization scheme, one must compute the true belief process \((\pi_t)_{t \ge 0}\) using the update equations, then quantize this to the set \(\mathcal{P}_N(\mathcal{X})\). This increases the computational complexity of this scheme (and requires explicit knowledge of the system model), both during the Q-learning algorithm and during implementation of the learned policy. Conversely, the sliding finite window scheme uses the approximate predictor from a fixed prior and a given history. Since there are only finitely many such histories, one can compute these offline before running the Q-learning and before implementation of the learned policy. They can then be accessed in a lookup table fashion.
    \item \emph{Model knowledge and a data-driven approach:} In the belief quantization scheme, both the encoder and the decoder must track the true belief \(\pi_t\). The encoder needs \(\pi_t\) to compute the proper value of \(\hat{\pi}_t\) and apply the learned policy, while the decoder needs it to compute the optimal reproduction value \(\hat{X}_t\). Thus, the belief quantization approach requires knowledge, or at least a good estimate, of the underlying model (in particular, the initial distribution and the kernels $T$ and $O$). Conversely, in the sliding finite window scheme the encoder policy is a constant function in \(\mathcal{P}(\mathcal{X})\), so it can directly use the sliding finite window to apply the learned policy. For the decoder, if one has knowledge of the model, it can simply take the form of \eqref{eq:decoder}, as we did in our implementation. However, theoretically one could also apply learning at the decoder: treat the decoder as a map $\gamma^d : W_t, Q_t, M'_t \to \hat{\mathcal{X}}$ and add this as another dimension in the Q-table, so that one now learns a joint encoder-decoder policy. Under the same Dobrushin coefficient conditions as for the encoder, this will become near-optimal for large \(N\). Of course, this comes at a significant memory and runtime cost for the Q-learning algorithm, so if one has model knowledge it is still beneficial to incorporate this at the decoder.
    \item \emph{Rate of convergence to near-optimality:} For the sliding finite window result in Theorem~\ref{theorem:finite window_near_optimal}, the convergence is exponential (note that although we only bound the expectation, since there are only finitely many initial memories, the convergence is also exponential for each initial finite window). Conversely, the belief quantization result in Corollary~\ref{corollary:quantization_near_optimal} is only asymptotic.
    \item \emph{On quantization:} For the belief quantization result in Corollary~\ref{corollary:quantization_near_optimal}, the quantization of \(\mathcal{P}(\mathcal{X})\) does not have to be uniform. Theoretically, this allows a more efficient quantization, although in our implementation we always use \cite[Algorithm 1]{Reznik} (which gives a uniform quantization). The sliding finite window can be seen as a non-uniform quantization of \(\mathcal{W}\) (since it is constant over the belief coordinate). However, it is uniform over the product space \((\mathcal{M}')^N \times \mathcal{Q}^N\), since the sliding finite window scheme uses every possible finite window.
\end{itemize}
We summarize some of the key differences between the schemes in Table \ref{tableComparisionRepeat}.


\begin{table}
\small
\begin{center}
    \begin{tabularx}{\textwidth} {
            >{\centering\arraybackslash}X
            || >{\centering\arraybackslash}X
            | >{\centering\arraybackslash}X
            | >{\centering\arraybackslash}X
            | >{\centering\arraybackslash}X
            | >{\centering\arraybackslash}X
            | >{\centering\arraybackslash}X}
                            & Convergence of Q-learning and near-optimality & Stability conditions not needed & Insensitive to initialization & Exponential convergence of performance & Bayesian update not needed & Lookup table implementation\\
        \hline
        Belief Quantization & \cmark                    & \cmark                 & \xmark                        & \xmark                      & \xmark       & \cmark                  \\
        \hline
        Sliding Finite Window       & \cmark                    & \xmark                 & \cmark                        & \cmark                      & \cmark & \cmark \end{tabularx}
\end{center}
\caption{Comparison of the two approximation schemes}\label{tableComparisionRepeat}
\vspace{-0.1in}
\end{table}

\normalsize

\subsection{Implications for the Average Cost Problem}
Under mild assumptions, it can be shown that as \(\beta \to 1\), a near-optimal policy for the discounted cost problem becomes near-optimal for the average cost problem (see \cite[Theorem 7.3.6]{yuksel2023control}). It is straightforward to show that these assumptions hold for the zero-delay coding MDP; indeed, this same set of assumptions was shown to hold in~\cite{wood2016optimal} and used to show the existence of optimal Walrand-Varaiya codes for the average cost problem. Then through~\cite[Theorem 7.3.6]{yuksel2023control} and our previous results we obtain the following corollaries. Note that the dependence on $\beta$ is only asymptotic; in practice we set $\beta = 0.9999$ and this seems to give policies very close to the optimum. For a fixed~$\beta$, the dependence on $\epsilon$ is as given in Theorems \ref{theorem:Qlearning2} and \ref{theorem:Qlearning1} (exponential for sliding finite window, asymptotic for belief quantization).
\begin{corollary}\label{corollary:avg2}
    Let \(\epsilon > 0\), let \(\widetilde{\gamma}^*_N\) be the policy given in Theorem~\ref{theorem:Qlearning2}, and assume the conditions of Theorem~\ref{theorem:Qlearning2} are met. Then there exists \(\beta \in (0,1)\) and \(N(\epsilon, \beta)\) sufficiently large such that for any \(\gamma\) which acts on the first \(N\) time steps to generate \(W_0\) and any \(\pi_{-N}\),
    \begin{align}
        \mathbf{E}^\gamma_{\pi_{-N}} \left[ J(W_0, \widetilde{\gamma}^*_N) - J^*(W_0)\right] < \epsilon.
    \end{align}
\end{corollary}

\begin{corollary}\label{corollary:avg1}
    Let \(\epsilon > 0\), let \(\widetilde{\gamma}^*_N\) be the policy given in Theorem~\ref{theorem:Qlearning1}. Then there exists \(\beta \in (0,1)\) and \(N(\epsilon, \beta)\) sufficiently large such that
    \begin{align}
        J(\zeta, \widetilde{\gamma}^*_N) - J^*(\zeta) < \epsilon.
    \end{align}
\end{corollary}

\section{Simulations}\label{section:simulations}
We now give some examples of zero-delay coding problems and simulate the performance of the policies resulting from Theorems~\ref{theorem:Qlearning1} and~\ref{theorem:Qlearning2}. In all of the following, we use a discount factor \(\beta = 0.9999\) and the distortion function \(d(x, \hat{x}) = (x - \hat{x})^2\). The average distortion is calculated over \(t=0\) to \(t=10^5\). The initial distribution is \(\zeta\) (that is, the invariant distribution for the source), so that the belief quantization scheme can be used. We observe that in each case with a known optimum, our Q-learning algorithm indeed reached this optimum (or very close to it), validating our theoretical results. In the case with an unknown optimum, we observe convergence rates which agree with the bounds in the theorems.

\emph{Remark:} For the implementation of the Q-learning algorithms, the theoretical upper bound for the possible number of states in each scheme grows very quickly in their respective parameters. In particular, for the belief quantization approach we have $\left| \mathcal{P}_N(\mathcal{X}) \right| = {N+|\mathcal{X}|+1 \choose |\mathcal{X}|-1}$ \cite{Reznik}, and for the sliding finite window approach we have $\left| 
\mathcal{W}_N \right| = \left| \mathcal{M}' \times \mathcal{Q}\right|^N$. However, the sets actually visited during Q-learning may be much smaller. Thus, one may wish to add entries to $V_t$ as they are visited by the Q-learning algorithm in a dynamic fashion. Furthermore, note that for certain problems it may be possible to significantly shrink the set of quantizers $\mathcal{Q}$ with no loss of optimality. For example, for a noiseless channel one can omit those quantizers with empty bins, or for an i.i.d. source those with non-convex bins. Due to this, the schemes are somewhat difficult to compare on the same set of axes, since the $N$ variable in each scheme leads to very different codebook sizes. Thus, when plotting against $N$ on the x-axis we provide the results on separate sets of axes.

Also note that for a stopping criteria in the algorithms, any measure of the convergence of \(V_t\) would be suitable; in our implementation, we stop when the relative pointwise difference is sufficiently small. For example, \(\max_{\hat{\pi},Q} (( V_{t+1}(\hat{\pi},Q) - V_t(\hat{\pi},Q)) / V_t(\hat{\pi},Q) < \epsilon\) for some small \(\epsilon > 0\).

\subsection{Comparison with Lloyd-Max quantizer for i.i.d.\ source and noiseless channel}
Let \(\mathcal{X} = \{0, \ldots, 7\}\) and consider an i.i.d.\ source \((X_t)_{t \ge 0}\), such that for all \(x\),
\begin{align}
    T(\cdot | x) = \begin{pmatrix}
                       \frac{1}{4} & \frac{1}{8} & \frac{1}{8} & \frac{1}{16} & \frac{1}{16} & \frac{1}{16} & \frac{1}{4} & \frac{1}{16}
                   \end{pmatrix}.
\end{align}
Note that in the i.i.d.\ case, we trivially have \(\delta(T) = 1\), so Corollary~\ref{corollary:dobrushin_finite_window} is applicable. Indeed, here we have that \(\alpha = 0\), so that for all \(N \ge 1\),
\begin{align}
    \mathbf{E}^{\gamma}_{\pi_{-N}} \left[ \left| J_\beta(W_0, \widetilde{\gamma}^*) - J^*_\beta(W_0) \right| \right] = 0.
\end{align}
That is, the optimal policy for the finite window representation is optimal (not just near-optimal) for the original problem for any \(N\). This is not surprising given the i.i.d.\ nature of the source; the approximation of \(\pi_{t-N}\) to \(\zeta\) is without any loss since \(\pi_{t-N}\) can be immediately recovered.

\begin{figure}[h]
    \centering
    \includegraphics[scale=0.85]{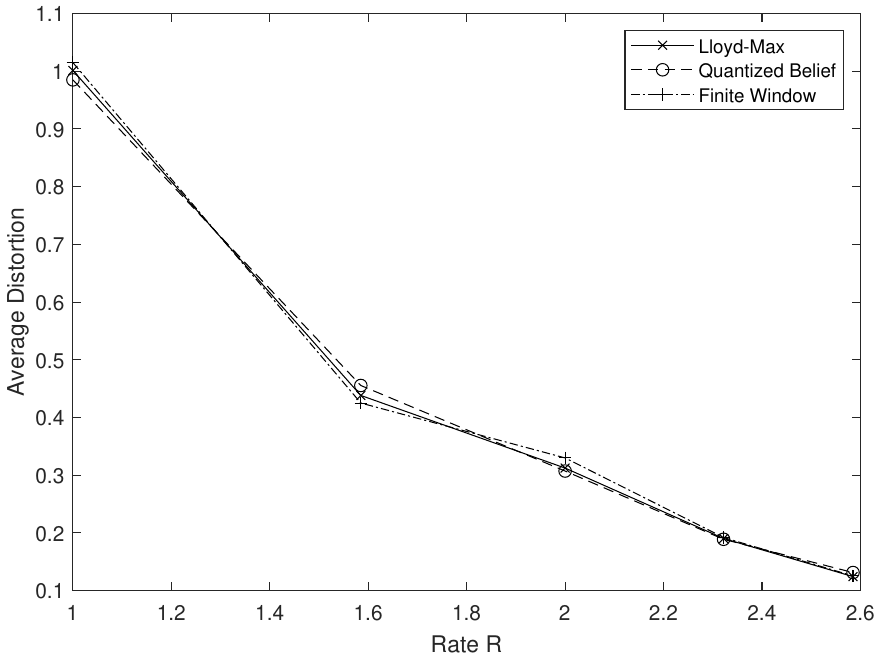}
    \caption{Comparison with Lloyd-Max}\label{fig:iid}
\end{figure}

Similarly for the quantization approach, \(N=1\) is sufficient since \(\pi_t = \zeta\) for all \(t \ge 0\), so increasing \(N\) does not change the resulting policy. Accordingly, we let \(N=1\) and compare the performance of both approaches against a Lloyd-Max quantizer when the channel is noiseless. We plot the performance for \(N=1\) and several sizes of \(\mathcal{M}\). The rate is calculated as
$R \coloneqq \log_2 |\mathcal{M}|$.
As expected by the above discussion, our algorithm matches with this quantizer in each case, which can be seen in Figure~\ref{fig:iid}.

\subsection{Comparison with memoryless encoding}
\begin{figure}[H]
    \centering
    \includegraphics[scale=0.85]{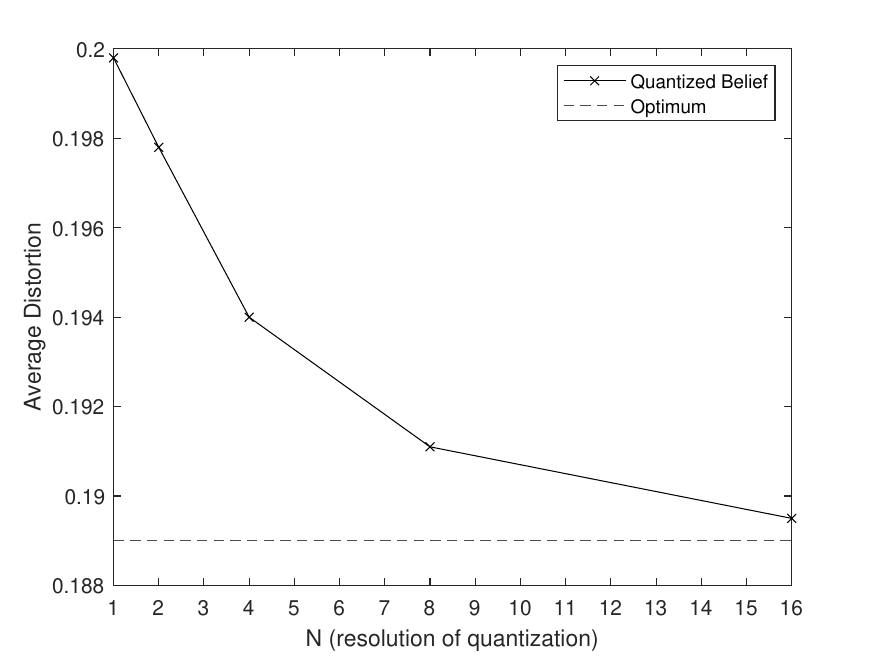}
    \caption{Quantized belief scheme vs memoryless encoding}\label{fig:quantized_memoryless}
\end{figure}

\begin{figure}[H]
    \centering
    \includegraphics[scale=0.85]{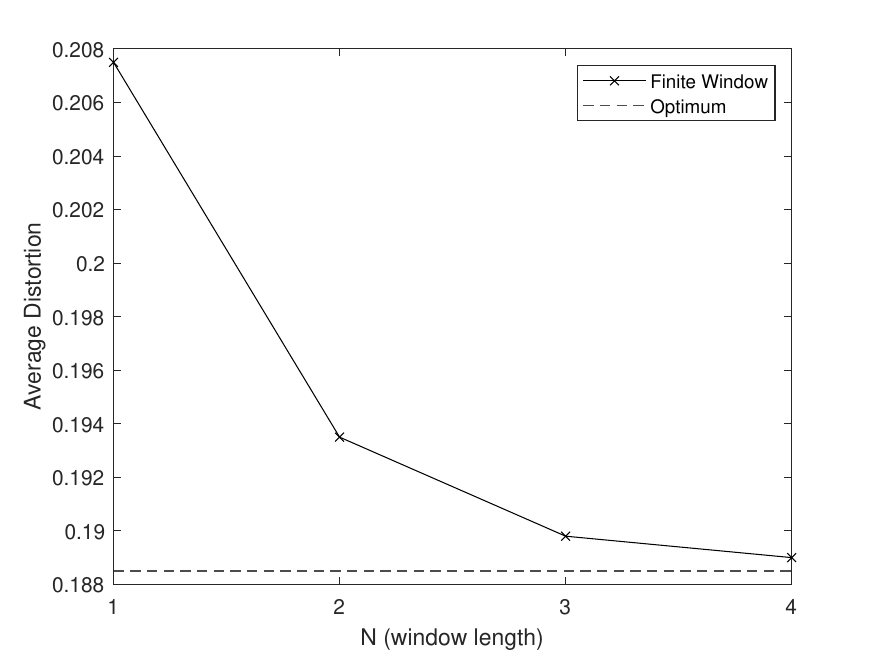}
    \caption{Finite memory scheme vs memoryless encoding}\label{fig:finite_memoryless}
\end{figure}

Consider now a Markov source with transition kernel given by
\begin{align}
    T = \begin{pmatrix}
            \frac{1}{2} & \frac{1}{6} & \frac{1}{6} & \frac{1}{6} \\
            \frac{1}{3} & \frac{1}{3} & \frac{1}{3} & 0           \\
            \frac{1}{3} & \frac{1}{3} & 0           & \frac{1}{3} \\
            \frac{1}{4} & \frac{1}{4} & \frac{1}{4} & \frac{1}{4}
        \end{pmatrix}
\end{align}
and where the channel is a 4-ary symmetric channel with error probability \(0.06\):
\begin{align}
    O = \begin{pmatrix}
            0.94 & 0.02 & 0.02 & 0.02 \\
            0.02 & 0.94 & 0.02 & 0.02 \\
            0.02 & 0.02 & 0.94 & 0.02 \\
            0.02 & 0.02 & 0.02 & 0.94
        \end{pmatrix}.
\end{align}
We have that \(\delta(T) = \frac{2}{3} > \frac{1}{2}\), so we can apply Corollary~\ref{corollary:dobrushin_finite_window}. In such a setup (where \(\mathcal{X} = \mathcal{M}\) and the channel is symmetric), it was shown in~\cite{WalrandVaraiya} that ``memoryless'' encoding (that is, where \(M_t = X_t\)) is optimal. We compare our algorithms against such an encoding policy, shown in Figures~\ref{fig:quantized_memoryless} and \ref{fig:finite_memoryless}, and note that both approach the optimum as \(N\) increases. Recall that $N$ represents {\em different} parameters in the different approximation schemes; for the quantized belief method, it represents the common denominator of the finite set $\mathcal{P}_N(\mathcal{X})$, while it represents the window length for the finite window method. Accordingly, we present the schemes on different sets of axes. Note that the convergence of the finite window scheme to the optimum indeed appears exponential in $N$, as expected by Theorem~\ref{theorem:Qlearning2}.

\subsection{Problem with unknown optimum}
Finally, we consider a setup where an optimal encoding scheme is unknown. Here we have a Markov source with transition kernel given by
\begin{align}
    T = \begin{pmatrix}
            0.2476 &   0.1527  &  0.0775  &  0.2219  &  0.2082 &   0.0920 \\
    0.0805  &  0.0247  &  0.0776 &   0.1290 &   0.3718 &   0.3164 \\
    0.1510  &  0.2335  &  0.2042 &   0.0107  &  0.1425 &   0.2580 \\
    0.0056  &  0.2252  &  0.2303  &  0.2173   & 0.1141   & 0.2076 \\
    0.1357  &  0.2685 &  0.0494   & 0.1981  &  0.2930  &  0.0553 \\
    0.2373  &  0.2795  &  0.0698  &  0.0399  &  0.1371  &  0.2363
        \end{pmatrix}.
\end{align}

\begin{figure}[h]
    \centering
    \includegraphics[scale=0.8]{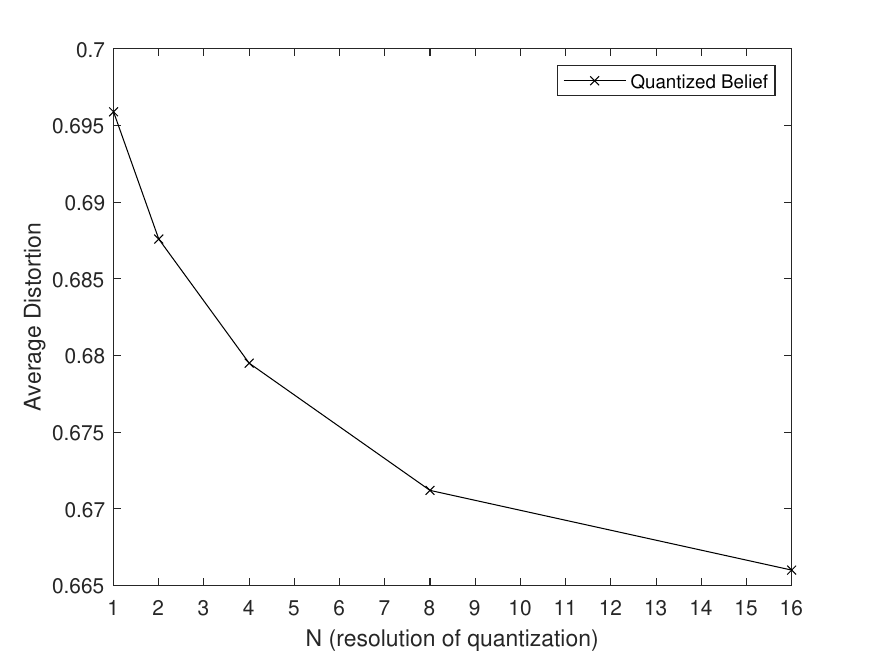}
    \caption{Quantized belief scheme with unknown optimum}\label{fig:quantized_unknown}
\end{figure}

\begin{figure}[h]
    \centering
    \includegraphics[scale=0.8]{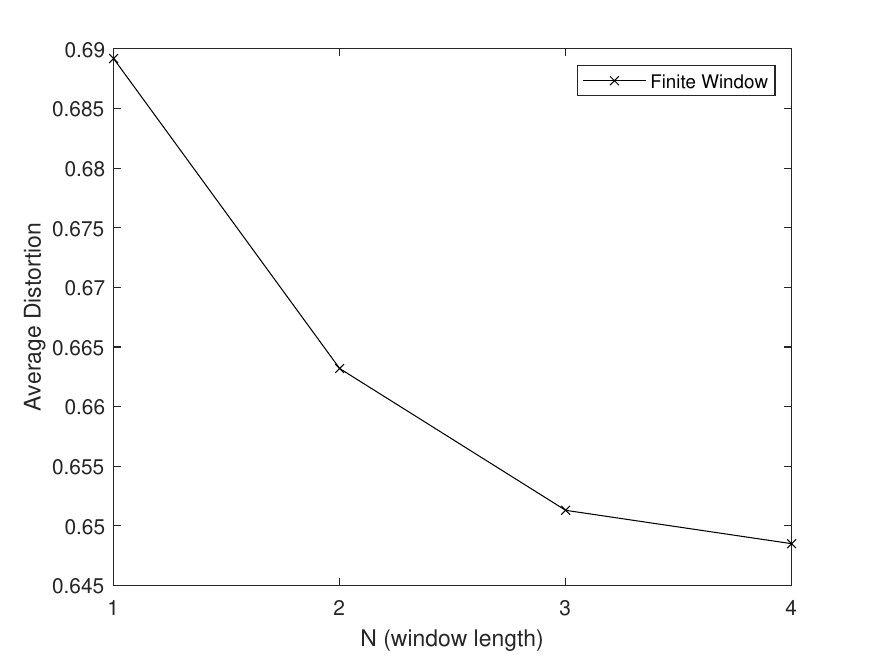}
    \caption{Finite memory scheme with unknown optimum}\label{fig:finite_unknown}
\end{figure}

Note that this transition kernel was randomly generated from the set of \(6 \times 6\) transition matrices. The channel is a 3-ary symmetric channel with error probability \(0.04\):
\begin{align}
    O = \begin{pmatrix}
            0.96 & 0.02 & 0.02 \\
            0.02 & 0.96 & 0.02 \\
            0.02 & 0.02 & 0.96
        \end{pmatrix}.
\end{align}
It can be verified that the Dobrushin coefficient conditions of Theorem~\ref{theorem:Qlearning2} are met for this setup. Figures~\ref{fig:quantized_unknown} and \ref{fig:finite_unknown} give the performance of both schemes for this problem.

\section{Conclusion}\label{section:conclusion}
We have provided two complementary approximation schemes for the zero-delay coding problem over a noisy channel with feedback, one in which the underlying belief space is quantized through a nearest-neighbor map, and the other in which the belief is approximated with a finite window of past observations. We then showed the convergence of a Q-learning algorithm to policies which are optimal for these approximations, yielding near-optimal policies for the zero-delay coding problem. Finally, we illustrated the convergence of our algorithm to the optimum through simulations. For future research, we wish to generalize our results to the case where our source-channel pair is continuous. We note that such a generalization should be mostly straightforward - under some recurrence assumptions on the source, one can still show the required regularity properties of the underlying MDP \cite{ghomi2021zero}, and then most of the approximation and convergence results go through with minor technical changes. We also wish to find less stringent filter stability conditions for the finite window scheme than the Dobrushin coefficient ones given here; this may be possible via observability-type conditions, such as those in~\cite{mcdonald2018stability}. The controlled case, where one communicates over a channel to a controller which affects some dynamical system, is also a promising research direction.

\appendices
\section{Proof of Theorems \ref{theorem:value_convergence} and \ref{theorem:finite window_near_optimal}}
\begin{lemma}\label{lemma:TV_bound}
    Consider the finite window approximation used in Section~\ref{section:finite window}, and recall \(I_t = (M'_{[t-N,t-1]}, Q_{[t-N,t-1]})\), \(W_t = (\pi_{t-N}, I_t)\), \(\hat{W}_t = (\zeta, I_t)\), and \(\hat{\pi}_t = P^\gamma_{\zeta}(X_t | m'_{[t-N,t-1]}, Q_{[t-N.t-1]})\). Then for any \(w_t \in \mathcal{W}\), \(\hat{w}_t \in \mathcal{W}_N\), and \(Q_t \in \mathcal{Q}\) we have
    \begin{equation}
        ||P(M'_t \in \cdot | w_t, Q_t) - P(M'_t \in \cdot | \hat{w}_t, Q_t)||_{\text{TV}} \le ||\pi_t - \hat{\pi}_t||_{\text{TV}}.
    \end{equation}
\end{lemma}
\begin{IEEEproof}
    Let \(f : \mathcal{M}' \to \mathbb{R}\) be measurable with \(||f||_\infty \le 1\). Then,
    \begin{align}
               & \left| \sum_{\mathcal{M}'}f(m'_t)P(m'_t | w_t, Q_t) - \sum_{\mathcal{M}'}f(m'_t)P(m'_t | \hat{w}_t, Q_t) \right| \\
        = \;   & \Biggl| \sum_{\mathcal{X}}\sum_{\mathcal{M}'}f(m'_t) P(m'_t | w_t, x_t, Q_t) P(x_t | w_t, Q_t) - \sum_{\mathcal{X}}\sum_{\mathcal{M}'}f(m'_t) P(m'_t | \hat{w}_t, x_t, Q_t) P(x_t | \hat{w}_t, Q_t) \Biggr|              \\
        = \;   & \Biggl| \sum_{\mathcal{X}}\sum_{\mathcal{M}'}f(m'_t) O_{Q_t}(m'_t | x_t) \pi_t(x_t)
        - \sum_{\mathcal{X}}\sum_{\mathcal{M}'}f(m'_t) O_{Q_t}(m'_t | x_t) \hat{\pi}_t(x_t)\Biggr|                                \\
        \le \; & ||\pi_t - \hat{\pi}_t||_{\text{TV}},
    \end{align}
    where the third line follows from conditional independence of \(M'_t\) and \(W_t\) given \((X_t, Q_t)\), and since \(Q_t\) is a function of \(W_t\) under any \(\gamma \in \Gamma_{\text{WS}}\). The last line follows from the fact that \(g(x) \coloneqq \sum_{\mathcal{M}'}f(m')O_Q(m'|x)\) is upper bounded by 1. Taking the supremum over all such \(f\) yields the result.
\end{IEEEproof}
\subsection*{Proof of Theorem \ref{theorem:value_convergence}}
We provide a proof for the case when \(N=1\), but an analogous proof follows for \(N>1\). It can be shown (see \cite[Theorem 4.2.3]{Hernandez}) that the functions \(J_\beta(w_t, \gamma)\) and \(J^*_\beta(w_t)\) satisfy the following fixed-point equations:
\begin{align}
    J_\beta(w_t, \gamma) & = c(w_t, \gamma(w_t)) + \beta \int_{\mathcal{W}} J_\beta(w_{t+1}, \gamma) P(dw_{t+1} | w_t, \gamma(w_t))                     \\
    J^*_\beta(w_t)       & = \min_{Q_t \in \mathcal{Q}} \Bigl( c(w_t, Q_t) + \beta \int_{\mathcal{W}} J^*_\beta(w_{t+1}) P(dw_{t+1} | w_t, Q_t) \Bigr),
\end{align}
for all \(w_t \in \mathcal{W}\) and \(\gamma \in \Gamma_{FS}\). Note that although the integral is over \(\mathcal{W}\), which is uncountable, we can only reach finitely many elements from a given \(w_t\) since the observation space \(\mathcal{M'}\) is finite. In particular, when \(N=1\), we can write \(w_t = (\pi_{t-1}, m'_{t-1}, Q_{t-1})\) and \(w_{t+1} = (\pi_t, m'_t, Q_t)\), so the above becomes
\begin{equation}
    J_\beta(w_t, \gamma) = c(w_t, \gamma(w_t)) + \beta \sum_{m'_t \in \mathcal{M'}} J_\beta((\pi_t, m'_t, \gamma(w_t)), \gamma)P(m'_t | w_t, \gamma(w_t)) \label{eq:policy-fixed-point}
\end{equation}
\begin{equation}
    J^*_\beta(w_t) = \min_{Q_t \in \mathcal{Q}} \Bigl( c(w_t, Q_t) + \beta \sum_{m'_t \in \mathcal{M'}} J^*_\beta(\pi_t, m'_t, Q_t) P(m'_t | w_t, Q_t) \Bigr).\label{eq:fixed-point}
\end{equation}

The functions \(\hat{J}_\beta(\hat{w}_t, \hat{\gamma})\) and \(\hat{J}^*_\beta(\hat{w}_t) \) satisfy equivalent fixed-point equations to~\eqref{eq:fixed-point}, so that for \(N=1\),
\begin{equation}
    \hat{J}_\beta(\hat{w}_t, \hat{\gamma}) = c_1(\hat{w}_t, \hat{\gamma}(\hat{w}_t)) + \beta \sum_{m'_t \in \mathcal{M'}} \hat{J}_\beta(\zeta, m'_t, \hat{\gamma}(\hat{w}_t))P(m'_t | \hat{w}_t, \hat{\gamma}(\hat{w}_t))
\end{equation}
\begin{equation}
    \hat{J}^*_\beta(\hat{w}_t) = \min_{Q_t \in \mathcal{Q}} \Bigl( c_1(\hat{w}_t, Q_t) + \beta \sum_{m'_t \in \mathcal{M'}} \hat{J}^*_\beta(\zeta, m'_t, Q_t) P(m'_t | \hat{w}_t, Q_t) \Bigr).\label{eq:fixed-point-approx}
\end{equation}

By definition of the extension \(\widetilde{J}_\beta^*\) we have \(\hat{J}^*_\beta(\hat{w}_1) = \widetilde{J}^*_\beta(w_1)\). Thus,
\begin{align}
    \beta \sum_{m'_0 \in \mathcal{M'}}\hat{J}^*_\beta(\zeta, m'_0, Q_0)P(m'_0| w_0, Q_0)
    = \beta \sum_{m'_0 \in \mathcal{M'}}\widetilde{J}^*_\beta(\pi_0, m'_0, Q_0)P(m'_0| w_0, Q_0).
\end{align}
We add and subtract the above term and use \(\widetilde{J}^*_\beta(w_0) = \hat{J}^*_\beta(\hat{w}_0)\) to obtain
\begin{align}
     & \left| \widetilde{J}^*_\beta(w_0) - J^*_\beta(w_0) \right| \\ = \; & \Bigl| \hat{J}^*_\beta(\hat{w}_0) - \beta \sum_{\mathclap{{m'_0 \in \mathcal{M'}}}}\hat{J}^*_\beta(\zeta, m'_0, Q_0)P(m'_0| w_0, Q_0) \\ & \qquad + \beta \sum_{\mathclap{{m'_0 \in \mathcal{M'}}}}\widetilde{J}^*_\beta(\pi_0, m'_0, Q_0)P(m'_0| w_0, Q_0) - J^*_\beta(w_0) \Bigr|.
\end{align}
Then applying the fixed-point equations~\eqref{eq:fixed-point} and~\eqref{eq:fixed-point-approx} to the last and first terms respectively,
\begin{multline*}
     \left| \widetilde{J}^*_\beta(w_0) - J^*_\beta(w_0) \right| \\ \hspace{-7cm} \le \max_{Q_0 \in \mathcal{Q}} \left| c_1(\hat{w}_0, Q_0) - c(w_0, Q_0)\right|  \\ + \beta \max_{Q_0 \in \mathcal{Q}} \left| \smashoperator[lr]{\sum_{m'_0}} \hat{J}^*_\beta(\zeta, m'_0, Q_0) P(m'_0 | \hat{w}_0, Q_0) - \smashoperator[lr]{\sum_{m'_0}}\hat{J}^*_\beta(\zeta, m'_0, Q_0)P(m'_0| w_0, Q_0) \right| \\ + \beta \max_{Q_0 \in \mathcal{Q}}\left| \smashoperator[lr]{\sum_{m'_0}}\widetilde{J}^*_\beta(\pi_0, m'_0, Q_0)P(m'_0| w_0, Q_0) - \smashoperator[lr]{\sum_{m'_0}} J^*_\beta(\pi_0, m'_0, Q_0) P(m'_0 | w_0, Q_0) \right|.
\end{multline*}
We now bound these three terms in expectation. The expectation is on \(W_0\) and \(\hat{W}_0\), with respect to the prior \(\pi_{-1}\) and some \(\gamma \in \Gamma_{\text{WS}}\), but we omit these in the expectation for notational simplicity. For the first term, by definition of \(c\) and \(c_N\) we have
\begin{equation}
    \mathbf{E}\left[ \left| c_1(\hat{W}_0, Q_0) - c(W_0, Q_0)\right| \right] \le ||d||_\infty \mathbf{E}\left[ ||\hat{\pi}_0 - \pi_0||_{\text{TV}} \right] \le ||d||_\infty L^1_0, \label{eq:cost_bound}
\end{equation}
where \(||d||_\infty = \max_{x, \hat{x}} d(x, \hat{x})\) and we recall the definition of \(L^N_t\) in~\eqref{eq:loss}; that is, \(L^1_0 = \sup_{\gamma \in \Gamma_{\text{WS}}}\mathbf{E}\left[ ||\pi_t - \hat{\pi}_t||_{\text{TV}} \right]\). For the second term, we have
\begin{align}
           \mathbf{E}&\left[ \max_{Q_0 \in \mathcal{Q}} \left| \smashoperator[lr]{\sum_{m'_0}} \hat{J}^*_\beta(\zeta, m'_0, Q_0) P(m'_0 | \hat{W}_0, Q_0) - \smashoperator[lr]{\sum_{m'_0}}\hat{J}^*_\beta(\zeta, m'_0, Q_0)P(m'_0| W_0, Q_0) \right| \right] \\
    \le \; & ||\hat{J}^*_\beta||_\infty \mathbf{E}\left[ \max_{Q_0 \in \mathcal{Q}} ||P(m'_0 | \hat{W}_0, Q_0) - P(m'_0| W_0, Q_0)||_{\text{TV}} \right]                                                                                                      \\
    \le \; & ||\hat{J}^*_\beta||_\infty L^1_0,
\end{align}
where the first inequality follows from the definition of the total variation (since \(\hat{J}^*_\beta / ||\hat{J}^*_\beta||_\infty\) is bounded by 1) and the second inequality is due to Lemma~\ref{lemma:TV_bound}. Finally, since both sums in the last term are over the same measure \(P(m'_0 | W_0, Q_0)\), we have
\begin{align}
           \mathbf{E}&\left[ \max_{Q_0 \in \mathcal{Q}}\left| \smashoperator[lr]{\sum_{m'_0}}\widetilde{J}^*_\beta(\pi_0, m'_0, Q_0)P(m'_0| W_0, Q_0) - \smashoperator[lr]{\sum_{m'_0}} J^*_\beta(\pi_0, m'_0, Q_0) P(m'_0 | W_0, Q_0) \right| \right] \\
    \le \; & \sup_{\gamma' \in \Gamma_{\text{WS}}} \mathbf{E}\left[ \left| \widetilde{J}^*_\beta(W_1) - J^*_\beta(W_1) \right| \right],
\end{align}
where we have used \((\pi_0, M'_0, Q_0) = W_1\). Combining all three bounds (and multiplying by \(\beta\) where appropriate) gives us
\begin{align}
    & \mathbf{E}\left[ \left| \widetilde{J}^*_\beta(W_0) - J^*_\beta(W_0) \right| \right] \\ \le & (||d||_\infty + \beta ||\hat{J}^*_\beta||_\infty)L^1_0 + \beta \sup_{\gamma' \in \Gamma_{\text{WS}}} \mathbf{E}\left[ \left| \widetilde{J}^*_\beta(W_1) - J^*_\beta(W_1) \right| \right]
\end{align}

We apply the same process to the final term, then recursively and by the fact that \(||J^*_\beta||_\infty \le \frac{||d||_\infty}{1-\beta}\), we have
\begin{equation}
    \mathbf{E} \left[ \left| \widetilde{J}_\beta^*(W_0) - J^*_\beta(W_0) \right| \right] \le \frac{||d||_\infty}{1-\beta}\sum_{t=0}^\infty \beta^t L^1_t.
\end{equation}
\hfill \(\blacksquare\)

\subsection*{Proof of Theorem \ref{theorem:finite window_near_optimal}}
As before, we consider \(N=1\), but analogous arguments follow for \(N>1\). We apply a similar strategy, by using the fixed-point equations in the proof of Theorem 2. Also, let \(Q_0^* \coloneqq \widetilde{\gamma}_1^*(w_0)\); that is, the action given by making the optimal policy for \(\text{MDP}_1\) constant over \(\mathcal{P}(\mathcal{X})\). Then, by~\eqref{eq:policy-fixed-point},
\begin{equation}
    J_\beta(w_0, \widetilde{\gamma}_1^*) = c(w_0, Q^*_0) + \beta \sum_{m'_0 \in \mathcal{M'}}J_\beta((\pi_0, m'_0, Q^*_0), \widetilde{\gamma}_1^*)P(m'_0 | w_0, Q^*_0)
\end{equation}
and using~\eqref{eq:fixed-point-approx} and the fact that \(\widetilde{J}^*_\beta(w_0) = \hat{J}^*_\beta(\hat{w_0})\),
\begin{align}
    \widetilde{J}^*_\beta(w_0) = c_1(\hat{w}_0, Q^*_0) + \beta \sum_{m'_0 \in \mathcal{M'}} \widetilde{J}^*_\beta(\pi_0, m'_0, Q_0^*) P(m'_0 | \hat{w}_0, Q_0^*).
\end{align}
Using \(w_1 = (\pi_0, m'_0, Q_0^*)\), we add and subtract
\begin{equation}
    \sum_{m'_0 \in \mathcal{M'}}\widetilde{J}^*_\beta(w_1)P(m'_0| w_0, Q^*_0),
\end{equation}
to obtain
\begin{align}
           & \left| J_\beta(w_0, \widetilde{\gamma}_1^*) - \widetilde{J}^*_\beta(w_0) \right| \\
    \le \; & \left| c(w_0, Q^*_0) - c(\hat{w}_0, Q^*_0) \right|
    \\ & + \beta \sum_{m'_0 \in \mathcal{M'}} \left| \widetilde{J}^*_\beta(w_1) P(m'_0 | \hat{w}_0, Q_0^*) - \widetilde{J}^*_\beta(w_1) P(m'_0 | w_0, Q_0^*) \right| \\ & + \beta \sum_{m'_0 \in \mathcal{M'}} \left| \widetilde{J}^*_\beta(w_1) - J_\beta(w_1, \widetilde{\gamma}_1^*) \right| P(m'_0 | w_0, Q^*_0).
\end{align}
Thus, using~\eqref{eq:cost_bound} and Lemma~\ref{lemma:TV_bound},
\begin{align}
        & \mathbf{E} \left[ \left| J_\beta(W_0, \widetilde{\gamma}_1^*) - \widetilde{J}^*_\beta(W_0) \right| \right]                                                                                     \\ \le \; & ||d||_\infty L^1_0                                                   \\
        & + \beta ||\widetilde{J}^*_\beta||_\infty \sup_{\gamma \in \Gamma_{\text{WS}}}\mathbf{E} \left[ ||P(m'_0 | \hat{W}_0, Q_0^*) - P(m'_0 | W_0, Q_0^*)||_{\text{TV}} \right]                                     \\
        & + \beta \sup_{\gamma \in \Gamma_{\text{WS}}}\mathbf{E} \left[ \left| \widetilde{J}^*_\beta(W_1) - J_\beta(W_1, \widetilde{\gamma}_1^*) \right| \right]                                                \\
    \le & (||d||_\infty + \beta ||\widetilde{J}^*_\beta||_\infty)L^1_0 + \beta \sup_{\gamma'}\mathbf{E} \left[ \left| J_\beta(W_1, \widetilde{\gamma}_1^*) - \widetilde{J}^*_\beta(W_1) \right| \right].
\end{align}
Recursively, and using the fact that \(||\widetilde{J}^*_\beta||_\infty \le \frac{||d||_\infty}{1-\beta}\),
\begin{align}
    \mathbf{E} \left[ \left| J^*_\beta(W_0, \widetilde{\gamma}_1^*) - J^*_\beta(W_0) \right| \right] & \le \frac{||d||_\infty}{1-\beta} \sum_{t=0}^\infty \beta^t L^1_t.\label{eq:policy_bound2}
\end{align}
Finally, we have
\begin{align}
           & \mathbf{E} \left[ \left| J^*_\beta(W_0, \widetilde{\gamma}_1^*) - J^*_\beta(W_0) \right| \right] \\ \le \; & \mathbf{E} \left[ \left| J^*_\beta(W_0, \widetilde{\gamma}_1^*) - \widetilde{J}^*_\beta(W_0) \right| \right] + \mathbf{E} \left[ \left| \widetilde{J}^*_\beta(W_0) - J^*_\beta(W_0) \right| \right] \\
    \le \; & \frac{2||d||_\infty}{1-\beta} \sum_{t=0}^\infty \beta^t L^1_t,
\end{align}
where the final inequality follows from~\eqref{eq:policy_bound2} and Theorem~\ref{theorem:value_convergence}.
\hfill \(\blacksquare\)

\section{Proof of Theorem \ref{theorem:dobrushin}}
\begin{lemma}\label{lemma:prior_bound}
    The following holds:
    \begin{equation}
        \sum_{\mathcal{X}} \sum_{\mathcal{M}'} ||\overline{\pi}^\mu_t - \overline{\pi}^\nu_t||_{\text{TV}} O_{Q_t}(m'|x) \pi^\mu_t(x) \le (2 - \tilde{\delta}(O))||\pi^\mu_t - \pi^\nu_t||_{\text{TV}},
    \end{equation}
    where \(\tilde{\delta}(O) = \min_{Q \in \mathcal{Q}}(\delta(O_Q))\).
\end{lemma}
\begin{IEEEproof}
     The following argument is from~\cite[Lemma 3.5]{mcdonald2020exponential}, adapted to our setup. Let \(f : \mathcal{X} \to \mathbb{R}\) be measurable with \(||f||_\infty \le 1\). Recall the update equations for \(\pi_t, \overline{\pi_t}\) given in~\eqref{eq:update}, and let $N^\mu(M_t,Q_t)$ denote the normalizing term for the $(\pi^\mu_t)_{t \ge 0}$ process, $N^\mu(M'_t,Q_t) \coloneqq \sum_{\mathcal{X}}g_{Q_t}(x, M'_t)\pi^\mu_t(x)$. Then we have for any \(M'_t = m'\) and \(Q_t = Q\),
    \begin{align}
               & \left| \sum_{\mathcal{X}}f(x)\overline{\pi}^\mu_t(x) - \sum_{\mathcal{X}}f(x)\overline{\pi}^\nu_t(x) \right|                                                                      \\
        =      & \left| \sum_{\mathcal{X}}f(x) \frac{g_{Q}(x, m')\pi^\mu_t(x)}{N^\mu(m', Q)} - \sum_{\mathcal{X}}f(x)\frac{g_{Q}(x, m')\pi^\nu_t(x)}{N^\nu(m', Q)}\right|          \\
        \le \; & \left| \sum_{\mathcal{X}}f(x) \frac{g_{Q}(x, m')\pi^\mu_t(x)}{N^\mu(m', Q)} - \sum_{\mathcal{X}}f(x)\frac{g_{Q}(x, m')\pi^\nu_t(x)}{N^\mu(m', Q)}\right|          \\
               & \hspace{3cm} + \left| \sum_{\mathcal{X}}f(x) \frac{g_{Q}(x, m')\pi^\nu_t(x)}{N^\mu(m', Q)} - \sum_{\mathcal{X}}f(x)\frac{g_{Q}(x, m')\pi^\nu_t(x)}{N^\nu(m', Q)}\right|        \\
        =  \;  & \frac{1}{N^\mu(m', Q)}
        \Biggl| \sum_{\mathcal{X}}f(x)g_{Q}(x, m')\pi^\mu_t(x) - \sum_{\mathcal{X}}f(x)g_{Q}(x, m')\pi^\nu_t(x) \Biggr|                                                                    \\
               & \hspace{3cm} + \Biggl| \frac{1}{N^\mu(m', Q)} - \frac{1}{N^\nu(m', Q)} \Biggr| \cdot \left|\sum_{\mathcal{X}}f(x)g_{Q}(x, m')\pi^\nu_t(x) \right|                                  \\
        \le    & \frac{1}{N^\mu(m', Q)} \sum_{\mathcal{X}}g_{Q}(x, m')|\pi^\mu_t - \pi^\nu_t|(x) \\ & \hspace{3cm} + \left| \frac{1}{N^\mu(m', Q)} - \frac{1}{N^\nu(m', Q)} \right| N^\nu(m', Q) \\
        \le    & \frac{1}{N^\mu(m', Q)} \Biggl( \sum_{\mathcal{X}}g_{Q}(x, m')|\pi^\mu_t - \pi^\nu_t|(x) + \left| N^\mu(m', Q) - N^\nu(m', Q) \right| \Biggr),
    \end{align}
    where in the second last line we have used the notation \(\sum_{\mathcal{X}} |\pi^\mu_t - \pi^\nu_t|(x) = \sum_{\mathcal{X}}(\mathbf{1}_{S^+} - \mathbf{1}_{S^-})(\pi^\mu_t - \pi^\nu_t)(x)\) with \(S^+ = \left\{ x |(\pi^\mu_t - \pi^\nu_t)(x) > 0 \right\}\) and \(S^- = \left\{ x | (\pi^\mu_t - \pi^\nu_t)(x) \le 0 \right\}\). Note that \(\sum_{\mathcal{X}} |\pi^\mu_t - \pi^\nu_t|(x) = ||\pi^\mu_t - \pi^\nu_t||_{\text{TV}}\). Taking the supremum over all \(f\) gives
    \begin{equation}
        ||\overline{\pi}^\mu_t - \overline{\pi}^\nu_t||_{\text{TV}} \le \frac{1}{N^\mu(m', Q)} \Biggl( \sum_{\mathcal{X}}g_{Q}(x, m')|\pi^\mu_t - \pi^\nu_t|(x) + \left| N^\mu(m', Q) - N^\nu(m', Q) \right| \Biggr). \label{eq:filter_bound}
    \end{equation}
    Thus, recalling that \(\psi\) is some appropriate reference measure over \(\mathcal{M}'\),
    \begin{align}
            & \sum_{\mathcal{X}} \sum_{\mathcal{M}'} ||\overline{\pi}^\mu_t - \overline{\pi}^\nu_t||_{\text{TV}} O_{Q_t}(m' | x)\pi^\mu_t(x)                                                                                            \\
        =   & \sum_{\mathcal{X}} \sum_{\mathcal{M}'} ||\overline{\pi}^\mu_t - \overline{\pi}^\nu_t||_{\text{TV}} g_{Q_t}(x, m')\psi(m')\pi^\mu_t(x)                                                                                     \\
        =   & \sum_{\mathcal{M}'} ||\overline{\pi}^\mu_t - \overline{\pi}^\nu_t||_{\text{TV}} \left( \sum_{\mathcal{X}} g_{Q_t}(x, m') \pi^\mu_t(x) \right) \psi(m')                                                                    \\
        =   & \sum_{\mathcal{M}'} ||\overline{\pi}^\mu_t - \overline{\pi}^\nu_t||_{\text{TV}} N^\mu(m', Q_t) \psi(m')                                                                                                                   \\
        \le & \sum_{\mathcal{M}'} \Bigl( \sum_{\mathcal{X}}g_{Q_t}(x, m')|\pi^\mu_t - \pi^\nu_t|(x) + \left| N^\mu(m', Q_t) - N^\nu(m', Q_t) \right| \Bigr)\psi(m')                                                            \\
        \le & \sum_{\mathcal{X}} \left( \sum_{\mathcal{M}'} g_{Q_t}(x, m') \psi(m') \right)|\pi^\mu_t - \pi^\nu_t|(x) \\ & \hspace{4cm} + \sum_{\mathcal{M}'} \left| \sum_{\mathcal{X}} g_{Q_t}(x, m')(\pi^\mu_t - \pi^\nu_t)(x)\right| \psi(m') \\
        \le & ||\pi^\mu_t - \pi^\nu_t||_{\text{TV}} + \sum_{\mathcal{M}'} \left| O_{Q_t}(\pi^\mu_t) - O_{Q_t}(\pi^\nu_t) \right|(m')                                                                                                  \\
        =   & ||\pi^\mu_t - \pi^\nu_t||_{\text{TV}} + ||O_{Q_t}(\pi^\mu_t) - O_{Q_t}(\pi^\nu_t)||_{\text{TV}},
    \end{align}
    where in the second last line we have used the kernel operator notation \(O_Q(\pi)(m') = \sum_{\mathcal{X}}O_Q(m' | x)\pi(x)\). It is shown in~\cite{dobrushin1956central} that the Dobrushin coefficient acts as a contraction coefficient for kernel operators under total variation. In particular
    \begin{equation}
        ||O_{Q_t}(\pi^\mu_t) - O_{Q_t}(\pi^\nu_t)||_{\text{TV}} \le (1-\delta(O_{Q_t}))||\pi^\mu_t - \pi^\nu_t||_{\text{TV}}.\label{eq:kernel_contraction}
    \end{equation}
    Thus,
    \begin{align}
            & \sum_{\mathcal{X}} \sum_{\mathcal{M}'} ||\overline{\pi}^\mu_t - \overline{\pi}^\nu_t||_{\text{TV}} O_{Q_t}(m'|x) \pi^\mu_t(x) \\ \le & (2 - \delta(O_{Q_t}))||\pi^\mu_t - \pi^\nu_t||_{\text{TV}}    \\
        \le & (2 - \tilde{\delta}(O))||\pi^\mu_t - \pi^\nu_t||_{\text{TV}},
    \end{align}
    where \(\tilde{\delta}(O) = \min_{Q \in \mathcal{Q}}(\delta(O_Q))\).
\end{IEEEproof}

\subsection*{Proof of Theorem~\ref{theorem:dobrushin}}
Note that, in \(\mathbf{E}_{\mu}^\gamma\) expectations of \(\overline{\pi}^\mu_t\) and \(\overline{\pi}^\nu_t\), it is enough to take the expectation over only \(M'_{[0,t]}\), since under any \(\gamma \in \Gamma_{\text{WS}}\), \(Q_{[0,t]}\) are deterministic given \(\mu\) and \(M'_{[0,t-1]}\). Thus,
\begin{align}
        & \mathbf{E}_{\mu}^\gamma\left[ ||\overline{\pi}^\mu_t - \overline{\pi}^\nu_t||_{\text{TV}} \right]                                                                                                                                   \\
    =   & \sum_{(\mathcal{M}')^{t+1}} ||\overline{\pi}^\mu_t - \overline{\pi}^\nu_t||_{\text{TV}} P^\gamma_{\mu}(m'_{[0,t]})                                                                                                                  \\
    =   & \sum_{(\mathcal{M}')^t} \sum_{\mathcal{X}} \sum_{\mathcal{M}'} ||\overline{\pi}^\mu_t - \overline{\pi}^\nu_t||_{\text{TV}} P^\gamma_{\mu}(m'_t | x_t, m'_{[0,t-1]}) P^\gamma_{\mu}(x_t | m'_{[0,t-1]}) P^\gamma_{\mu}(m'_{[0,t-1]}) \\
    =   & \sum_{(\mathcal{M}')^t} \sum_{\mathcal{X}} \sum_{\mathcal{M}'} ||\overline{\pi}^\mu_t - \overline{\pi}^\nu_t||_{\text{TV}} O_{Q_t}(m'_t | x_t) \pi^\mu_t(x) P^\gamma_{\mu}(m'_{[0,t-1]})                                            \\
    \le & (2 - \tilde{\delta}(O)) \sum_{(\mathcal{M}')^t}||\pi^\mu_t - \pi^\nu_t||_{\text{TV}} P^\gamma_{\mu}(m'_{[0,t-1]})                                                                                                                   \\
    =   & (2 - \tilde{\delta}(O)) \mathbf{E}_{\mu}^\gamma \left[ ||\pi^\mu_t - \pi^\nu_t||_{\text{TV}} \right],
\end{align}
where the third equality follows from the fact that \(Q_t\) is deterministic given \(\mu\) and \(M'_{[0,t-1]}\), and that given \(X_t\) and \(Q_t\), \(M'_t\) depends only on the kernel \(O_{Q_t}\). For the inequality we used Lemma~\ref{lemma:prior_bound}. Finally, using the Dobrushin contraction property for kernels (as noted in the derivation of \eqref{eq:kernel_contraction}), we have
\begin{align}
    \mathbf{E}_{\mu}^\gamma \left[ ||\pi^\mu_{t+1} - \pi^\nu_{t+1}||_{\text{TV}} \right] & \le (1-\delta(T))\mathbf{E}_{\mu}^\gamma\left[ ||\overline{\pi}^\mu_t - \overline{\pi}^\mu_t||_{\text{TV}} \right]     \\
                                                                                  & \le (1-\delta(T))(2 - \tilde{\delta}(O)) \mathbf{E}_{\mu}^\gamma \left[ ||\pi^\mu_t - \pi^\nu_t||_{\text{TV}} \right].
\end{align}
\hfill \(\blacksquare\)

\section{Proof of Theorem \ref{theorem:Qlearning2}}
We will show certain ergodic behavior (in particular, \cite[Assumption 2.1]{kara2023qlearning}) and then invoke \cite[Theorem 2.1]{kara2023qlearning}. Indeed, we have trivially that
\begin{align}
    \frac{\sum_{k=0}^t c_N(\hat{w}_k, Q_k) \mathbf{1}(\hat{w}_k = \hat{w}, Q_k = Q)}{\sum_{k=0}^t\mathbf{1}(\hat{w}_k = \hat{w}, Q_k = Q)} = c_N(\hat{w},Q),
\end{align}
so that \cite[Assumption 2.2 (ii)]{kara2023qlearning} holds. Additionally, by positive Harris recurrence of \((X_t)_{t \ge 0}\), the marginals on \(P(X_t \in \cdot)\) converge to \(\zeta\), so we have for any \(f\) that
\begin{align}
    \frac{\sum_{k=0}^t f(\hat{w}_{k+1}) \mathbf{1}(\hat{w}_k = \hat{w}, Q_k = Q)}{\sum_{k=0}^t\mathbf{1}(\hat{w}_k = \hat{w}, Q_k = Q)} \to \int f(\hat{w}_1)P_N(\hat{w}_1 | \hat{w}, Q)
\end{align}
almost surely as \(t \to \infty\), where \(P_N\) is from~\eqref{eq:eta_hat}. Thus \cite[Assumption 2.1 (iii)]{kara2023qlearning} holds.

Finally, note that not every \(\hat{w} \in \mathcal{W}_N\) has positive probability of being visited (certain sequences of channel outputs and quantizers are impossible depending on the source and channel), so we restrict ourselves to only those with positive probability. Then \cite[Assumption 2.1 (i)]{kara2023qlearning} holds for (almost every) \(\hat{w}\), so we apply \cite[Theorem 2.1]{kara2023qlearning} to obtain that almost surely, as \(N \to \infty\),
\begin{equation}
    V_t(\hat{w},Q) \to V^*(\hat{w},Q),
\end{equation}
where
\begin{equation}
    V^*(\hat{w},Q) = c_N(\hat{w},Q) + \beta \sum_{\hat{w}_1 \in \mathcal{W}_N} \min_Q V^*(\hat{w}_1) P_N(\hat{w}_1 | \hat{w}, Q),
\end{equation}
where \(P_N\) and \(c_N\) are from \eqref{eq:eta_hat} and \eqref{eq:c_hat}. Taking the minimum over \(\mathcal{Q}\) gives us the classic Bellman optimality equation, and hence \(\hat{J}_\beta(\hat{w}, \hat{\gamma}^*_N) = \hat{J}^*_\beta(\hat{w})\) for almost every \(\hat{w} \in \mathcal{W}_N\). The result follows by applying Corollary \ref{corollary:dobrushin_finite_window}. \hfill \(\blacksquare\)

\bibliographystyle{IEEEtran}
\bibliography{CDC_ACC_Edited,SerdarBibliography}

\end{document}